\documentclass[12pt]{article}

\usepackage{amsthm,amssymb,amsmath}
\usepackage{graphicx}
\usepackage{mdframed}
\usepackage{fancybox}
\usepackage{xcolor}
\usepackage[pagebackref=true,colorlinks,linkcolor=blue,citecolor=magenta]{hyperref}
\usepackage{tocbibind}
\usepackage{makeidx}
\makeindex
\theoremstyle{definition} 
\newtheorem{definition}{Definition}[section]
\theoremstyle{theorem}
\newtheorem{theorem}[definition]{Theorem}
\newtheorem{lemma}[definition]{Lemma}
\newtheorem{proposition}[definition]{Proposition}
\newtheorem{corollary}[definition]{Corollary}
\newtheorem{remark}[definition]{Remark}
\theoremstyle{definition}

\begin{document}

\title{The general form of $p$-completely contractive homomorphisms of the $p$-analog of the Fourier-Stieltjes algebras}
\date{}
\maketitle

\begin{center}
\author{Mohammad Ali Ahmadpoor \\Department of Pure Mathematics, Faculty of Mathematical Sciences, University of Guilan, Rasht, Iran.\\  email \href{mailto: m-a-ahmadpoor@phd.guilan.ac.ir}{m-a-ahmadpoor@phd.guilan.ac.ir}\\
\and Marzieh Shams Yousefi \\Department of Pure Mathematics, Faculty of Mathematical Sciences, University of Guilan, Rasht, Iran.\\ email \href{mailto:m.shams@guilan.ac.ir}{m.shams@guilan.ac.ir} }
\end{center}

\begin{abstract}
In this paper, we follow two main goals. In the first attempt, we give some functorial properties of the $p$-analog of the Fourier-Stieltjes algebras in which we generalize some previously existed definitions and theorems in Arsac and Cowling's works, to utilize them to prove $p$-complete boundedness of some well-known maps on these algebras. In the second part, as an application of these generalizations, we prove $p$-completely boundedness of homomorphisms which are induced by continuous and proper piecewise affine maps that is a generalization of Ilie's work on Fig\`a-Talamanca-Herz algebras.
\end{abstract}
\textbf{keywords:} {Completely bounded homomorphisms, Fourier-Stieltjes algebras, $QSL_p$-spaces, Piecewise affine maps}        



\section{Introduction}\label{subsection1.1}

Let $G$ be a locally compact group. The Fourier algebra, $A(G)$, and the Fourier-Stieltjes algebra, $B(G)$, on the locally compact group $G$, have been found by Eymard in 1964 \cite{EYMARD1964}. The general form of special type of maps on the Fourier and Fourier-Stieltjes algebras has been studied extensively. For example, when $G$ is an Abelian topological group, $A(G)$ is nothing except $L_1(\widehat{G})$, where $\widehat{G}$ is the Pontrjagin dual group of $G$, and $B(G)$ is isometrically isomorphic to $M(\widehat{G})$, the measure algebra. In this case, Cohen in \cite{COHEN1960-1} and \cite{COHEN1960-2} studied homomorphisms from $L_1(G)$ to $M(H)$, for Abelian groups $G$ and $H$, and gave the general form of these maps, as the weighted maps by a piecewise affine map on the underlying groups.\\
By \cite{BLECHER1992,EFFROSRUAN1991}, we know that $A(G)$ and $B(G)$ are operator spaces as the predual of a von Neumann algebra, and the dual of a $C^*$-algebra, respectively. Ilie in \cite{ILIE2004} and \cite{ILIESPRONK2005} studied the completely bounded homomorphisms from the Fourier to the Fourier-Stieltjes algebras. It is shown that for a continuous piecewise affine map $\alpha:Y\subset H\rightarrow G$, the homomorphism $\Phi_\alpha:A(G)\rightarrow B(H)$, defined through
\begin{align*}
\Phi_\alpha u=\left\{
\begin{array}{ll}
u\circ\alpha& \text{on} \: Y\\
0& \text{o.w.}
\end{array}\right. ,\quad u\in A(G),
\end{align*}
is  completely bounded. Moreover, in the cases that $\alpha$ is an affine map and a homomorphism, the homomorphism $\Phi_\alpha$ is completely contractive and completely positive, respectively.

The Fig\`a-Talamanca-Herz algebras were introduced by Fig\`a-Talamanca for Abelian locally compact groups \cite{FIGATALAMANCA1965}, and it is generalized for any locally compact group by Herz \cite{HERZ1971}. For $p\in (1,\infty)$, coefficient functions of the left regular representation of a locally compact group $G$ on $L_{p}(G)$ give us the Fig\`a-Talamanca-Herz algebra $A_p(G)$, and we have $A_2(G)=A(G)$. Therefore,  Fig\`a-Talamanca-Herz algebras can be seen as the $p$-analog of the Fourier algebras.

Daws in \cite{DAWS2010} introduced the $p$-operator space structure, with an extensive application to $A_p(G)$, which generalizes the operator space structure of $A(G)$.\\
Oztop and Spronk in \cite{OZTOPSPRONK2012}, and Ilie in \cite{ILIE2013} studied the $p$-completely bounded homomorphisms on the Fig\`a-Talamanca-Herz algebras, using the $p$-operator space structure. In \cite{ILIE2013} it is shown that the map $\Phi_\alpha :A_p(G)\rightarrow A_p(H)$, defined via
\begin{align*}
\Phi_\alpha u= \left\{
\begin{array}{ll}
u\circ \alpha & \text{on} \: Y\\
0 & \text{o.w}
\end{array}\right.,\quad u\in A_p(G),
\end{align*}
is a $p$-completely (bounded) contractive homomorphism for a continuous proper (piecewise) affine map $\alpha : Y\subset H\rightarrow G $ in the case that the locally compact group $H$ is amenable.

Runde in \cite{RUNDE2005} found $p$-analog of the Fourier-Stieltjes algebras, $B_p(G)$. He used extensively the theory of $ QSL_p $-spaces and representations on these spaces. Also, he gave the $p$-operator space structure of $B_p(G)$. More precisely, it is indicated that the space $B_p(G)$ is the dual space of the algebra of universal $p$-pseudofunctions $UPF_p(G)$, i.e. $B_p(G)=UPF_p(G)^*$. Therefore, by combining this result with the fact that for a concrete operator space like $UPF_p(G)$, we have $B_p(G)=\mathcal{CB}_p(UPF_p(G),\mathbb{C})$ \cite{DAWS2010}, it can be concluded that $B_p(G)$ is a $p$-operator space as a dual of a $p$-operator space. The second author of this paper studied the $p$-analog of the Fourier-Stieltjes algebras on the inverse semigroups in \cite{shams}.

In this paper, for a continuous proper map $\alpha : Y\subset H\rightarrow G $, we study the weighted maps $\Phi_\alpha : B_p(G)\rightarrow B_p(H)$ which is defined by
\begin{align}\label{eq1}
\Phi_\alpha u= \left\{
\begin{array}{ll}
u\circ \alpha & \text{on} \: Y\\
0 & \text{o.w}
\end{array}\right.,\quad u\in B_p(G).
\end{align}
We will show that when $\alpha$ is an affine map, $\Phi_\alpha$ is a $p$-complete contraction, and in the case that $\alpha$ is a piecewise affine map, it is $p$-completely bounded homomorphism. For this aim, we put amenability assumption on open subgroups of $H$. Our approach to the concept of $p$-operator space structure on the $p$-analog of the Fourier-Stieltjes algebra, is the $p$-operator structure that can be implemented on this space from its predual.

The paper is organized as follows: First we give required  definitions and theorems about the $p$-analog of the Fourier-Stieltjes algebras and representations on $QSL_p$-spaces in Section \ref{SECTIONPRELIMINARIES}. In Section \ref{SECTIONGENERALIZATIONS}, first we generalize Cohen-Host type idempotent theorem on the Fourier-Stieltjes algebras (see \cite{RUNDE2007}) to the $p$-analog of the Fourier-Stieltjes algebras  in Theorem \ref{IDMTHEOREM}. In addition, by Theorem \ref{PROPOSITIONDUALOFPFP} we give general form of the main theorem in \cite{COWLING1979}. For this aim, we need to give $p$-analog of some results in \cite{ARSAC1976}. As a crucial theorem in this paper, we have Theorem \ref{IMPORTANTMAPPING} in Section \ref{APPLICABLEHOMO} which will be applied in the next section. Final section, Section \ref{WEIGHTEDHOMO}, is about generalization of Ilie's results on homomorphisms of the Fig\`a-Talamanca-Herz algebras in \cite{ILIE2013}. Indeed, we study $p$-completely boundedness of homomorphisms of the form \eqref{eq1}.

\section{Preliminaries}\label{SECTIONPRELIMINARIES}
In this paper, $G$ and $H$ are locally compact groups, and for $p\in (1,\infty)$, the number $p'$ is its complex conjugate, i.e. $1/p+1/p'=1$. In the first step, we give essential notions and definitions on $QSL_p$-spaces, and representations of groups on such spaces. For more information one can see \cite{RUNDE2005}.
\begin{definition}\label{def1}
A representation of a locally compact group $ G $ is a
pair $ (\pi , E) $, where $ E $ is a Banach space and $ \pi $ is a group homomorphism from $ G $ into the
invertible isometries on $ E $, that is continuous with respect to the given topology on $ G $
and the strong operator topology on $ \mathcal{B}(E) $.
\end{definition}

\begin{remark}
Every representation $ (\pi , E) $ of a locally compact group $ G $ induces a representation of the group algebra $L_1(G)$ on $ E $, i.e. a contractive algebra homomorphism from $L_1(G)$ into
$ \mathcal{B}(E) $, which we shall denote likewise by $ \pi $, through
\begin{align}\label{11}
&\pi(f)=\int f(x)\pi (x)dx,\ f\in L_1(G),\\
&\langle \pi(f)\xi, \eta\rangle =\int f(x)\langle \pi (x)\xi, \eta\rangle dx,\quad \xi\in E,\:\eta\in E^*,\nonumber
\end{align}
where the integral \eqref{11} converges with respect to the strong operator topology.
\end{remark}

\begin{definition}
Let $ (\pi , E)$ and $(\rho , F) $ be representations of the locally compact group $G$. Then

\begin{enumerate}
\item
$ (\pi , E)$ and $ (\rho , F) $ are called equivalent, if there exists an invertible isometry $\varphi : E\rightarrow F$ such that
\begin{equation*}
\varphi\pi (x)\varphi^{-1}=\rho (x),\qquad x\in G.
\end{equation*}
\item
$ (\rho , F) $ is said to be a subrepresentation of $ (\pi , E)$, if $F$ is a closed subspace of $E$, and for every $x\in G$ we have $\pi(x){|_{F}}=\rho(x)$.
\item
$ (\rho , F) $ is said to be contained in $ (\pi , E)$, if it is equivalent to a subrepresentation of $  (\pi , E) $, and will be denoted by $(\rho , F) \subset (\pi , E)$.
\end{enumerate}
\end{definition}

\begin{definition}
\begin{enumerate}
\item A Banach space is called an $L_p$-space if it is of the form $L_p(X)$ for some measure
space $X$.
\item A Banach space is called a $ QSL_p $-space if it is isometrically isomorphic to a quotient of a subspace of an $ L_p $-space.
\end{enumerate}
\end{definition}
We denote by $ \text{Rep}_p(G) $ the collection of all (equivalence classes) of representations of $G $ on a  $ QSL_p $-space.

\begin{definition}
A representation of a Banach algebra $\mathcal{A}$ is a pair $ (\pi , E)$, where $E$ is a Banach space, and $\pi $ is a contractive algebra homomorphism from $\mathcal{A}$ to $\mathcal{B}(E)$. We call
$ (\pi , E)$ isometric if $\pi $ is an isometry and essential if the linear span of $\{\pi(a)\xi \: : \ a \in \mathcal{A},\: \xi\in
E\} $ is dense in $E$. 
\end{definition}

\begin{remark}\label{rem11}
If $ G $ is a locally compact group and $ (\pi , E)$ is a representation of $ G $ in the sense of Definition
\ref{def1}, then \eqref{11} induces an essential representation of $L_1(G)$. Conversely,
every essential representation of $L_1(G)$ arises in this fashion.
\end{remark}

\begin{definition}

\begin{enumerate}
\item
A representation $ (\pi , E)\in \text{Rep}_p(G) $ is called cyclic, if there exists $\xi_0\in E$ such that $\pi (L_1(G))\xi_0$ is dense in $E$. The set of cyclic representations of group $G$ on $QSL_p$-spaces is denoted by $\text{Cyc}_p(G)$.
\item
A representation $ (\pi , E)\in \text{Rep}_p(G)$ is called $p$-universal, if it contains every cyclic representation.
\end{enumerate}

\end{definition}

\begin{remark}\label{REMARKGARDELLA}
We know  that every $p$-universal representation of $G$, contains every cyclic representation of $G$ on a $QSL_p$-space, in the sense of equivalency. In Addition, every representation in $\text{Rep}_p(G)$ is contained in a $p$-universal representation. Actually, one could make a new $p$-universal representation by constructing a direct sum of an arbitrary representation with a $p$-universal representation. For more about representations of a locally compact group $G$
on $L _p$-spaces and $QSL_p$-spaces see \cite{GARDELLATHIEL2015}.
\end{remark}

Now we are ready to describe the Fig\`a-Talamanca-Herz, and the $p$-analog of the Fourier-Stieltjes algebras.

\begin{definition}
Fig\`a-Talamanca-Herz algebra on the locally compact group $G$, which is denoted by $A_p(G)$, is the collection of functions $u:G\rightarrow\mathbb{C}$ of the form

\begin{equation}\label{Ap}
u(\cdot)=\sum_{n=1}^{\infty}\langle \lambda_{p}(\cdot)\xi_n,\eta_n\rangle,
\end{equation}
with 

\begin{equation}\label{Ap2}
(\xi_n)_{n\in\mathbb{N}}\subset L_{p}(G),\quad (\eta_n)_{n\in\mathbb{N}}\subset L_{p'}(G),\quad\text{and}\quad\sum_{n=1}^\infty\|\xi_n\| \|\eta_n\|<\infty,
\end{equation}
where $ \lambda_{p}$ is the left regular representation of $G$ on $L^{p}(G)$, defined as
\begin{align*}
&\lambda_{p}: G\rightarrow \mathcal{B}(L^{p}(G)),\quad\lambda_{p}(x)\xi(y)=\xi(x^{-1}y),\quad \xi\in L^{p}(G),\:x,y\in G.
\end{align*}
The norm of $A_p(G)$ is defined as 
\begin{equation*}
\| u\|=\inf\Big\{\sum_{n=1}^\infty\|\xi_n\| \|\eta_n\|\: :\ u(\cdot)=\sum_{n=1}^\infty\langle \lambda_{p}(\cdot)\xi_n, \eta_n\rangle\Big\},
\end{equation*}
where the infimum is taken over all expressions of $ u $ in \eqref{Ap} with \eqref{Ap2}.
With this norm and pointwise operations, $A_p(G)$ turns into a commutative regular Banach algebra.
\end{definition}

\begin{remark}
 The $p$-analog of the Fourier-Stieltjes algebra has been studied, for example in \cite{COWLING1979}, \cite{FORREST1994}, \cite{MIAO1996} and \cite{PIER1984}, as the multiplier algebra of the Fig\`a-Talamanca-Herz algebra. In this paper, we follow the construction of Runde in definition and notation (See \cite{RUNDE2005}) which we swap indexes $p$ and $p'$.
\end{remark}

\begin{definition}\label{definitionofBp}
The set of all functions of the form

\begin{equation*}
u(\cdot)=\langle \pi (\cdot)\xi, \eta\rangle,\ \xi\in E,\ \eta\in E^*,\quad\text{for some}\ (\pi , E)\in \text{Rep}_{p}(G),
\end{equation*}
equipped with the norm
\begin{equation*}
\| u\|=\inf\Big\{\sum_{n=1}^\infty\|\xi_n\| \|\eta_n\|\: :\ u(\cdot)=\sum_{n=1}^\infty\langle \pi_n (\cdot)\xi_n, \eta_n\rangle \Big\},
\end{equation*}
where 
\begin{align*}
(\pi_n,E_n)_{n\in\mathbb{N}}\subset\text{Cyc}_p(G),\quad\text{with}\quad \sum_{n=1}^\infty\|\xi_n\| \|\eta_n\|<\infty,
\end{align*}
 is denoted by $B_p(G)$, and is called the $p$-analog of the Fourier-Stieltjes algebra of the locally compact group $G$.
\end{definition}

\begin{remark}\label{REMARKRUNDERUNDE}
\begin{enumerate}
\item
By \cite[Lemma 4.6]{RUNDE2005}, the space $B_p(G)$ can be defined to be the set of all coefficient functions of a $p$-universal representation $(\pi,E)$, and the norm of an element $u\in B_p(G)$ is the infimum of all values $\sum_{n=1}^\infty\|\xi_n\|\|\eta_n\|<\infty$, which such vectors exist in the representation of $u$ as a coefficient function of $(\pi,E)$, i.e. $u(\cdot)=\sum_{n=1}^\infty\langle\pi(\cdot)\xi_n,\eta_n\rangle$.
\item

By \cite[Theorem 4.7]{RUNDE2005}, the space $B_p(G)$ equipped with the norm defined as above, and pointwise operations is a commutative unital Banach algebra, and by \cite[Corollary 5.3]{RUNDE2005}, by denoting multiplier algebra of $A_p(G)$ by $\mathcal{M}(A_p(G))$, we have the following contractive embeddings
\begin{align*}
A_p(G)\subset B_p(G)\subset \mathcal{M}(A_p(G)).
\end{align*}
\item\label{ITEMFIGARESTRICTION} In \cite[Lemma 2.4]{RUNDE2007}, the following identification is shown for an open subgroup $G_0$ of a locally compact group $G$  
\begin{align*}
A_p(G_0)\cong\{f\in A_p(G)\ :\ \text{supp}(f)\subset G_0\},
\end{align*}
and through this fact, one can assume that functions in $A_p(G_0)$ are restriction of functions in $A_p(G)$ to the open subgroup $G_0$.
\end{enumerate}
\end{remark}

\begin{definition}
Let $ (\pi , E)\in \text{Rep}_p(G) $.
\begin{enumerate}
\item For each $f\in L_1(G)$, let $\| f\|_{\pi}:=\|\pi(f)\|_{\mathcal{B}(E)}$, then $\|\cdot\|_{\pi}$ defines an algebra seminorm on $L_1(G)$.
\item By $ PF_{p,\pi}(G) $, we mean the $ p $-pseudofunctions associated with $ (\pi , E) $, which is the closure of
$ \pi( L_1(G)) $ in $ \mathcal{B}(E) $.
\item
If $ (\pi , E) = (\lambda_p,L_p(G))$, we denote $PF_{p,\lambda_p}(G)$ by  $PF_p(G)$.
\item
If $ (\pi , E) $ is $ p $-universal, we denote $ PF_{p,\pi}(G) $ by $ UPF_{p}(G)$, and call it the algebra of
universal $ p $-pseudofunctions.
\end{enumerate}
\end{definition}
\begin{remark}
\begin{enumerate}
\item
For $ p = 2 $, the algebra $ PF_{p}(G) $ is the reduced group $ C^* $-algebra, and $ UPF_{p}(G)$ is the
full group $ C^* $-algebra of $ G $.
\item
If $ (\rho , F)\in \text{Rep}_p(G) $ is such that $ (\pi , E) $ contains every cyclic subrepresentation of $ (\rho , F)$, then $\|\cdot\|_{\rho}\leq\|\cdot\|_{\pi}$ holds. In particular, the definition of $ UPF_{p}(G)$ is independent of a particular $ p $-universal representation.
\item
With $\langle\cdot,\cdot\rangle$ denoting $L_1(G)-L_\infty(G)$ duality, and with $(\pi,E)$ a $p$-universal representation of $G$, we have
\begin{align*}
\| f\|_\pi=\sup\{|\langle f,g\rangle|\ :\ g\in B_p(G),\ \|g\|_{B_p(G)}\leq 1\},\quad f\in L_1(G).
\end{align*}
\end{enumerate} 
\end{remark}

Next lemma states that $ B_p(G) $  is a dual space.
\begin{lemma}{\cite[Lemma 6.5]{RUNDE2005}}\label{RUNDEDUALITY}
Let $ (\pi , E)\in \text{Rep}_{p}(G)$. Then, for each $  \phi\in PF_{p,\pi}(G)^*$, there is a unique $g\in B_p(G) $, with $\|g\|_{B_p(G)}\leq \|\phi\| $ such that
\begin{equation}\label{duality}
\langle\pi (f),\phi\rangle=\int_G f(x)g(x)dx,\qquad f\in L_1(G).
\end{equation}
Moreover, if $ (\pi , E) $ is $p$-universal, we have $\|g\|_{B_p(G)}= \|\phi\| $.
\end{lemma}

The $p$-operator space structure which is used in this paper is Daws' approach for $A_p(G)$ \cite{DAWS2010}. A concrete $p$-operator space is a closed subspace of $\mathcal{B}(E)$, for some $QSL_p$-space $E$. In this case for each $n\in\mathbb{N}$ one can define a norm $\|\cdot\|_n$ on $\mathbb{M}_n(X)=\mathbb{M}_n\otimes X$ by identifying $ \mathbb{M}_n(X) $ with a subspace of $ \mathcal{B}(l_p^n\otimes_p E $). So, we have the family of norms $\Big(\|\cdot\|_n\Big)_{n\in\mathbb{N}}$ satisfying:

\begin{enumerate}
\item[$\mathcal{D}_\infty:$] For $u\in \mathbb{M}_n(X)$ and $v\in \mathbb{M}_m(X)$, we have that $ \|u\oplus v\|_{n+m}=\max\{\|u\|_n,\|v\|_m\} $. Here $u\oplus v\in \mathbb{M}_{n+m}(X)$ has block representation
$\begin{pmatrix}
u & 0 \\
0 & v
\end{pmatrix}. $

\item[$\mathcal{M}_p:$] For every $u\in \mathbb{M}_m(X)$ and $\alpha\in \mathbb{M}_{n,m}$, $\beta\in \mathbb{M}_{m,n}$ we have that 
$$ \|\alpha u\beta\|_n\leq\|\alpha\|_{\mathcal{B}(l^m_p,l^n_p)}\|u\|_m\|\beta\|_{\mathcal{B}(l^n_p,l^m_p)}. $$

\end{enumerate}

\begin{definition}
A linear operator $\Psi :X\rightarrow Y$ between two $p$-operator spaces is called $p$-completely bounded, if $\|\Psi\|_{\text{p-cb}}=\sup_{n\in\mathbb{N}}\|\Psi^{(n)}\|<\infty$, and $p$-completely contractive if $\|\Psi\|_{\text{p-cb}}=\sup_{n\in\mathbb{N}}\|\Psi^{(n)}\|\leq 1$, where $\Psi^{(n)} :\mathbb{M}_n(X)\rightarrow\mathbb{M}_n(Y)$ is defined in the natural way.
\end{definition}

\begin{theorem}\cite[Theorem 4.3]{DAWS2010}\label{DAWSTHEOREM}
Let $X$ be a $p$-operator space. There exists a $p$-complete isometry
$\varphi : X^*\rightarrow\mathcal{B}(l_p(I))$ for some index set $I$.
\end{theorem}

\begin{lemma}\cite[Lemma 4.5]{DAWS2010}\label{DAWSLEMMA}
If $ \Psi:X\rightarrow Y $ is $p$-completely bounded map between two operator spaces $X$ and $Y$, then $\Psi^* :Y^*\rightarrow X^*$ is $p$-completely bounded, with $\|\Psi^*\|_{\text{p-cb}}\leq \|\Psi\|_{\text{p-cb}}$.
\end{lemma}

\begin{remark}
\begin{enumerate}
\item
It should be noticed that, converse of Lemma \ref{DAWSLEMMA} is not necessarily true, unless $X$ be a closed subspace of $\mathcal{B}(E)$, for some $L_p$-space $E$.
\item
In the case of $p$-analog of the Fourier-Stieltjes algebras, from duality $B_p(G)=UPF_p(G)^*$, and the fact that $UPF_{p}(G)\subset\mathcal{B}(E)$, for $p$-universal representation $(\pi,E)$, by Theorem \ref{DAWSTHEOREM}, we can induce $p$-operator space structure from predual $UPF_p(G)$ on $B_p(G)$ through identification $B_p(G)=\mathcal{CB}_p(UPF_p(G),\mathbb{C})$.
\item In comparison to \cite{ILIE2013}, because of above explanations, a major difference in our work is that we need to study predual of some crucial $p$-completely bounded maps (See Theorem \ref{IMPORTANTMAPPING}), instead of their duals.
\end{enumerate}
\end{remark}

In Section \ref{WEIGHTEDHOMO}, we will study the homomorphisms on the $p$-analog of the Fourier-Stieltjes algebras induced by the continuous map $\alpha :Y\subset H\rightarrow G$, in the cases that $\alpha$ is homomorphism, affine and piecewise affine map, and $Y$ in the coset ring of $H$. So, we give some preliminaries here.\\
For a locally compact topological group $H$, let $\Omega_0(H)$ denote the ring of subsets which generated by open cosets of $H$. By \cite{ILIE2013} we have
\begin{align}\label{OPENCOSETRINGS}
\Omega_0(H)=\left\{Y\backslash\cup_{i=1}^nY_i \: :
\begin{array}{ll}
&Y\:\text{is an open coset of }\: H,\\
& Y_1,\ldots ,Y_n\:\text{open subcosets of infinite index in}\: Y\\
\end{array}\right\}.
\end{align}
Moreover, for a set $Y\subset H$, by $\text{Aff}(Y)$ we mean the smallest coset containing $Y$, and if $Y=Y_0\backslash\cup_{i=1}^nY_i\in\Omega_0(H)$, then $\text{Aff}(Y)=Y_0$.
Similarly, let us denote by  $\Omega_{\text{am-}0}(H)$ the ring of open cosets of open amenable subgroups of $H$, i.e.
\begin{align}\label{OPENCOSETRINGS}
\Omega_{\text{am-}0}(H)=\left\{Y\backslash\cup_{i=1}^nY_i \: :
\begin{array}{ll}
&Y\:\text{is an open coset of an open amenable subgroup of }\: H,\\
& Y_1,\ldots ,Y_n\:\text{open subcosets of infinite index in}\: Y\\
\end{array}\right\}.
\end{align}

\begin{definition}\label{DEFPIECEWISE}
Let $\alpha : Y\subset H\rightarrow G$ be a map.
\begin{enumerate}
\item
The map $\alpha$ is called an affine map on an open coset $Y$ of an open subgroup $H_0$, if
\begin{equation*}
\alpha(xy^{-1}z)=\alpha(x)\alpha(y)^{-1}\alpha(z),\qquad x,y,z\in Y,
\end{equation*}

\item
The map $\alpha$ is called a piecewise affine map if
\begin{enumerate}
\item
there are pairwise disjoint $ Y_i\in\Omega_0(H)$, for $ i=1,\ldots , n $, such that $Y=\cup_{i=1}^nY_i$,

\item there are affine maps $\displaystyle{\alpha_i : \text{Aff}(Y_i)\subset H\rightarrow G}$, for $ i=1,\ldots , n $, such that
\begin{equation*}
\alpha |_{Y_i}=\alpha_i |_{Y_i}.
\end{equation*}

\end{enumerate}
\end{enumerate}

\end{definition}

\begin{definition}
If $X$ and $Y$ are locally compact spaces, then a map $\alpha :Y\rightarrow X$ is called proper, if $\alpha^{-1}(K)$ is compact subset of $Y$, for every compact subset $K$ of $ X $.
\end{definition}

\begin{proposition}\cite[Proposition 4]{DUNKLRAMIREZ1971}\label{ramirez}
Let $\alpha : H\rightarrow G$ be a continuous group homomorphism. Then $\alpha$ is proper if and only if the bijective homomorphism $\tilde{\alpha}: H/{\ker\alpha}\rightarrow \alpha(H)=G_0$, is a topological group isomorphism, and $\ker\alpha$ is compact.
\end{proposition}

\begin{remark}\label{AFFFIINEREMMM}
\begin{enumerate}
\item
Proposition \ref{ramirez} implies that every continuous proper homomorphism is automatically a closed map. Therefore, $\alpha(H)$ is a closed subgroup of $G$. Also, $\ker\alpha$ is a compact normal subgroup of $H$.
\item
It is well-known that $\tilde{\alpha}$ is a group isomorphism, if and only if $\alpha$ is an open homomorphism into $\alpha(H)$, with the relative topology.  
\item\label{affine-remark}{\cite[Remark 2.2]{ILIE2004}}
If $Y=h_0H_0$ is an open coset of an open subgroup $H_0\subset H$, and $\alpha : Y\subset H\rightarrow G$ is an affine map, then there exists a group homomorphism $ \beta $ associated to $\alpha$ such that
\begin{align}\label{affine-homomorphism}
&\beta : H_0\subset H\rightarrow G,\quad\beta (h)= \alpha(h_0)^{-1}\alpha(h_0h),\quad h\in H_0.
\end{align}
\item \label{HOMAFFPROPER}
It is clear that, $\alpha$ is a proper affine map, if and only if $\beta$ is a proper homomorphism.
\item\cite[Lemma 8]{ILIE2013}\label{LEMMA8888} Let $Y\in\Omega_0(H)$, and $\alpha:\text{Aff}(Y)\rightarrow G$ be an affine map such that $\alpha|_{Y}$ is proper, then $\alpha$ is proper.
\end{enumerate}
\end{remark}

\section{Some generalizations}\label{SECTIONGENERALIZATIONS}

In order to being prepared for Theorem \ref{IDMTHEOREM} which is a generalization of \cite[Theorem 1.5]{RUNDE2007}, we need some elementary definitions and facts which we give them in the following.

\begin{definition}\label{DEFINITIONUNICONV}
\begin{enumerate}
\item
A Banach space $(E, \|\cdot\|)$ is said to be uniformly convex if for
every $0 < \epsilon\leq 2$ there is $\delta > 0$ so that for any two vectors $x$ and $y$ in $E$ with $\|x\|=\|y\| = 1$, the condition $\|x-y\| \geq\epsilon$ implies that $\|\frac{x+y}{2}\|\leq 1-\delta$. Intuitively, the center of a line segment inside the unit ball must lie deep inside
the unit ball unless the segment is short.
\item
A Banach space $E$ is said to be smooth if for each $\xi\in E\backslash\{0\}$ there exists a unique $\eta\in E^*$ such that $\|\eta\|=1$ and $\langle\xi,\eta\rangle =\|\xi\|$.

\end{enumerate}

\end{definition}
\begin{remark}\label{REMUNICONV}
It is worthwhile to note that by Definition \ref{DEFINITIONUNICONV}, every closed
subspace of a uniformly convex Banach space is again a uniformly convex Banach
space.
\end{remark}
Now we state an immensely important theorem about a quotient space which
can be found in \cite{ISTRATESCU1983}.
\begin{theorem}\cite[Theorem 2.4.18]{ISTRATESCU1983}\label{THEOREMQUOTUNI}
Let $E$ be a uniformly convex Banach space and $F$ be a closed linear subspace of $E$. Then the quotient space $E/F$ is uniformly convex Banach space.
\end{theorem}
Now we can conclude the following statement.
\begin{corollary}\label{QSLPSMOOUNICON}
 Every $QSL_p$-space $E$ is uniformly convex and smooth.
 
\begin{proof}
Uniformly convexity of $QSL_p$-space $E$ can be derived from Remark \ref{REMUNICONV}
and Theorem \ref{THEOREMQUOTUNI}. Since $E$ is uniformly convex, by \cite[Lemma 8.4(i) and
Theorem 9.10]{FHHMPZ2001} it is concluded that $E^*$ is smooth, but $E^*$ is a $QSL_{p'}$-space so is
uniformly convex, and then $E^{**}$ is smooth, but $E = E^{**}$ so $E$ is smooth.
\end{proof} 

\end{corollary}

\begin{theorem}\label{IDMTHEOREM}
For a subset $C\subset G$ following statements are equivalent.
\begin{enumerate}
\item\label{idm1}
$C$ is a left open coset,
\item\label{idm2}
$ \chi_C\in B(G) $ with $ \|\chi_C\|_{B(G)}=1 $,
\item\label{idm3}
$\chi_C\neq 0$ is a normalized coefficient function of a representation $(\pi,E)$ where $E$ or $E^*$ is smooth,
\item\label{idm4}
$ \chi_C\in B_p(G) $ with $ \|\chi_C\|_{B_p(G)}=1 $.
\end{enumerate}
\begin{proof}
Equivalency of the first three statements have been proved in \cite[Theorem 1.5]{RUNDE2007}. We demonstrate \eqref{idm2}$\Rightarrow$\eqref{idm4}$\Rightarrow$\eqref{idm3}. Let \eqref{idm2} hold. Then from the fact that $B(G)\subset B_p(G)$ and this embedding is a contraction, we have $\chi_C\in B_p(G)$ with $ \|\chi_C\|_{B_p(G)}\leq 1$, which by inequality $\|\cdot\|_{C_b(G)}\leq \|\cdot\|_{B_p(G)}$, we have $ \|\chi_C\|_{B_p(G)}=1$ which shows \eqref{idm2} implies \eqref{idm4}.\\
Now let $\chi_C\in B_p(G)$ with $ \|\chi_C\|_{B_p(G)}= 1$. So, by Definition \ref{definitionofBp}, the function $\chi_C$ is a normalized coefficient function of an isometric group representation on a $QSL_{p}$-space, which is smooth by Corollary \ref{QSLPSMOOUNICON} that is \eqref{idm3}.
\end{proof}

\end{theorem}

\begin{corollary}\label{COROIDEMCOSETRING}
Let $G$ be a locally compact group and $Y\in\Omega_0(G)$, then we have $\chi_Y\in B_p(G)$. Moreover, we have
\begin{align}\label{FFF}
1\leq\|\chi_Y\|_{B_p(G)}\leq 2^{m_Y},\quad\text{with}\quad m_Y=\inf\{m\in\mathbb{N}\; :\ Y=Y_0\backslash\cup_{i=1}^mY_i\},
\end{align}
where for $i=0,1,\ldots,m$ sets $Y_i$, are as \eqref{OPENCOSETRINGS}.
\begin{proof}
Since $Y\in\Omega_0(G)$, then by \eqref{OPENCOSETRINGS}, there exist open coset $Y_0$ and open subcosets $Y_i\subset Y_0$, for $i=1,\ldots , m$ and $m\in\mathbb{N}$ such that $Y=Y_0\backslash\cup_{i=1}^mY_i$. By Theorem \ref{IDMTHEOREM}-\eqref{idm4}, we have $\chi_{Y_i}\in B_p(G)$, with $\|\chi_{Y_i}\|_{B_p(G)}=1$, for $i=0,1,\ldots ,m$. On the other hand, since 
\begin{align}\label{EQCHI}
\chi_Y=\chi_{Y_0}-\big(\sum_{i=1}^m\chi_{Y_i}-\sum_{i,j=1}\chi_{Y_i\cap Y_j}+\sum_{i,j,k=1}\chi_{Y_i\cap Y_j\cap Y_k}+\cdots+(-1)^{m+1}\chi_{Y_1\cap Y_2\cap\ldots\cap Y_m}\big),
\end{align}
then we have $\|\chi_Y\|_{B_p(G)}\leq 2^m$, and by taking infimum on all possible decomposition of $Y$ as \eqref{OPENCOSETRINGS} relation \eqref{FFF} holds.
\end{proof}

\end{corollary}

In the sequel, we will give some extensions of results in \cite{ARSAC1976}. For a unitary representation $(\pi,\mathcal{H}_\pi)$ with Hilbert space $\mathcal{H}_\pi$, the $\pi$-Fourier space  has been defined to be closed linear span of the set of the coefficient functions of the representation $(\pi,\mathcal{H}_\pi)$, and is denoted by $A_\pi$ with the norm in usual way. Moreover, $\pi$-Fourier-Stieltjes algebra, $B_\pi$, for such representation is defined to be $w^*$-closure of $A_\pi$. Additionally, if we let $C^*_\pi(G)$ be the  $C^*$-algebra associated with $\pi$, we have $B_\pi=C^*_\pi(G)^*$. Here we introduce $p$-generalization of these results.

\begin{definition}
For a representation $(\pi , E)\in\text{Rep}_p(G)$, we define the \textit{$p$-analog of the $\pi$-Fourier space}, $A_{p,\pi}$, to be closed linear span of the collection of the coefficient functions of representation $(\pi, E)$ equipped with the norm 
\begin{align*}
\| u\|_{A_{p,\pi}}=\inf\{\sum_{n=1}^\infty\|\xi_n\|\|\eta_n\|\ :\  u(\cdot)=\langle\pi(\cdot)\xi_n,\eta_n\rangle,\ (\xi_n)_{n\in\mathbb{N}}\subset E, \; (\eta_n)_{n\in\mathbb{N}}\subset E^*\},\quad u\in A_{p,\pi}.
\end{align*}
\end{definition}

\begin{remark}\label{REMARKPSIMAP}
\begin{enumerate}
\item\label{REMARKPSIMAP1}
Let $(\pi , E)\in\text{Rep}_p(G)$. Consider the following map,
\begin{align*}
&\Psi_{p,\pi}: E^*\widehat{\otimes}E\rightarrow C_b(G),\qquad \Psi_{p,\pi}\bigg(\sum_{n=1}^\infty\eta_n\otimes\xi_n\bigg)=\sum_{n=1}^\infty\langle\pi(\cdot)\xi_n,\eta_n\rangle,
\end{align*}
so we can identify coimage of $\Psi_{p,\pi}$ with Banach space $E^*\widehat{\otimes}E/\ker{\Psi_{p,\pi}}$, which implies that the norm on coimage is the quotient norm i.e.,
\begin{align*}
\| \sum_{n=1}^\infty\xi_n\otimes\eta_n+\ker\Psi_{p,\pi}\|&=\inf\{\sum_{n=1}^\infty\|x_n \| \| y_n \| \ : \sum_{n=1}^\infty\langle\pi(\cdot)x_n,y_n\rangle=\sum_{n=1}^\infty\langle\pi(\cdot)\xi_n,\eta_n\rangle \}\\
&=\|\sum_{n=1}^\infty\langle\pi(\cdot)\xi_n,\eta_n\rangle \|_{A_{p,\pi}}.
\end{align*}
So, one can identify $A_{p,\pi}$ with coimage of $\Psi_{p,\pi}$, or equivalently with the quotient space $E^*\widehat{\otimes}E/\ker\Psi_{p,\pi}$.
\item
Since we have $A_{p,\pi}\cong E^*\widehat{\otimes}E/\ker\Psi_{p,\pi}$, then the space $A_{p,\pi}$ is a Banach space.

\end{enumerate}
\end{remark}

\begin{proposition}\label{LEMMANORMAPP}

Let $(\pi , E)\in\text{Rep}_p(G)$ and $\text{Cyc}_{p,\pi}(G)=\{ (\rho , F)\in\text{Cyc}_{p}(G)\ \&\ (\rho , F)\subset (\pi , E)\}$. Then for $u(\cdot)=\sum_{n=1}^\infty\langle\pi(\cdot)\xi_n,\eta_n\rangle\in A_{p,\pi}$ we have
\begin{align*}
\|u\|_{A_{p,\pi}}=\inf\{\sum_{n=1}^\infty\|x_n \| \| y_n \| \ : \ u(\cdot)=\sum_{n=1}^\infty\langle\rho_n(\cdot)x_n,y_n\rangle  \},
\end{align*}
where the infimum is taken on all representations of $u$ in which $((\rho_n, F_n))_{n\in\mathbb{N}}\subset\text{Cyc}_{p,\pi}(G)$ with $(x_n)_{n\in\mathbb{N}}\subset F_n $ and $(y_n)_{n\in\mathbb{N}}\subset F^*_n$.
\begin{proof}

Let
\begin{equation*}
C:=\inf\{\sum_{n=1}^\infty\|x_n\|\|y_n\|\; :\ u(\cdot)=\sum_{n=1}^\infty\langle\rho_n(\cdot)x_n,y_n\rangle,\ (\rho_n,F_n)_{n\in\mathbb{N}}\subset \text{Cyc}_{p,\pi}(G)\}.
\end{equation*}
Assume that $u(\cdot)=\sum_{n=1}^\infty\langle\pi(\cdot)\xi_n,\eta_n\rangle$ with $\sum_{n=1}^\infty\|\xi_n\|\|\eta_n\|<\infty$. For each $n\in\mathbb{N}$ we may put
\begin{align*}
&F_n=\overline{\pi(L_1(G))\xi_n{}{}}^{\|\cdot\|_{E}},\quad \rho_n :G\rightarrow\mathcal{B}(F_n),\quad \rho_n(x)=\pi(x)|_{F_n},\quad
\quad x_n=\xi_n,\ y_n=\eta_n|_{F_n},
\end{align*}
then we have
\begin{align*}
&((\rho_n, F_n))_{n\in\mathbb{N}}\subset\text{Cyc}_{p,\pi}(G),\quad u(\cdot)=\sum_{n=1}^\infty\langle\rho_n(\cdot)x_n,y_n\rangle,
\end{align*}
with $ C\leq\sum_{n=1}^\infty\|x_n\|\|y_n\|\leq \sum_{n=1}^\infty\|\xi_n\|\|\eta_n\|$.
Since $(\xi)_{n\in\mathbb{N}}\subset E$ and $(\eta)_{n\in\mathbb{N}}\subset E^*$ are arbitrary in the representing of $u$, we have $C\leq\| u\|_{A_{p,\pi}}$.

For the inverse inequality, let $\epsilon >0$ is given. Then there exist $((\rho_n, F_n))_{n\in\mathbb{N}}\subset\text{Cyc}_{p,\pi}(G)$, $(x_n)_{n\in\mathbb{N}}\subset F_n$, $(y_n)_{n\in\mathbb{N}}\subset F_n^*$, and for each $n\in\mathbb{N}$, we have $(\rho_n, F_n)\subset (\pi,E)$ such that
\begin{equation*}
\sum_{n=1}^\infty\|x_n\|\|y_n\| < C+\epsilon,\qquad u(\cdot)=\sum_{n=1}^\infty\langle\rho_n(\cdot)x_n,y_n\rangle .
\end{equation*}
Now for each $n\in\mathbb{N}$, by applying Hahn-Banach theorem extend each $y_n\in F_n^*$ to the $\eta_n\in E^*$ such that $\|\eta_n\|=\|y_n\|$. Therefore,
\begin{equation*}
\| u\|_{A_{p,\pi}}\leq \sum_{n=1}^\infty\|x_n\|\|\eta_n\| =\sum_{n=1}^\infty\|x_n\|\|y_n\| < C+\epsilon,
\end{equation*}
and it means $\| u\|_{A_{p,\pi}}\leq C$.

\end{proof}

\end{proposition}

For a representation $(\pi,E)\in \text{Rep}_p(G)$, by $\big(\pi^\infty, l_p(\mathbb{N},E)\big)$ we denote the representation
\begin{align*}
\pi^\infty:G\rightarrow\mathcal{B}(l_p(\mathbb{N},E)),\quad\pi^\infty(x)((\xi_n)_{n=1}^\infty)=(\pi(x)\xi_n)_{n=1}^\infty,\ x\in G,\ (\xi_n)_{n=1}^\infty\in l_p(\mathbb{N},E).
\end{align*}
Similarly, for a free ultrafilter $\mathcal{U}$, and ultrapower of the space $l_p(\mathbb{N},E)$, by $\big((\pi^\infty)_\mathcal{U}, l_p(\mathbb{N},E\big)_\mathcal{U})$, we mean the representation
\begin{align*}
(\pi^\infty)_\mathcal{U}: G\rightarrow\mathcal{B}(l_p(\mathbb{N},E)_\mathcal{U}),\quad(\pi^\infty)_\mathcal{U}(x)((\xi_n)_\mathcal{U})=(\pi(x)\xi_n)_\mathcal{U},\ x\in G,\ (\xi_n)_\mathcal{U}\in l_p(\mathbb{N},E)_\mathcal{U}.
\end{align*}

\begin{proposition}\label{PROPULTRAAP}
For each $(\pi, E)\in\text{Rep}_p(G)$, there exists a free ultrafilter $\mathcal{U}$, such that by restricting $(\pi^\infty)_{\mathcal{U}}$ to the subspace 
\begin{align*}
F=\overline{\{(\pi^\infty)_{\mathcal{U}}(f)(x)\ :\ f\in L_1(G)\ x\in l_p(\mathbb{N},E)_\mathcal{U}\}{}}^{\|\cdot\|_{l_p(\mathbb{N},E)_{\mathcal{U}}}}\subset l_p(\mathbb{N},E)_{\mathcal{U}},
\end{align*}
the representation $((\pi^\infty)_\mathcal{U},F)$ is weak-weak$^*$ continuous, essential and isometric representation of $PF_{p,\pi}(G)$, and we have $PF_{p,\pi}(G)^*=\overline{A_{p,(\pi^\infty)_\mathcal{U}}{}}^{w^*}$.

\begin{proof}

By \cite[Lemma 6.5]{RUNDE2005}, there exists a free ultrafilter $\mathcal{U}$ such that canonical representation of $PF_{p,\pi}(G)$ on $l_p(\mathbb{N},E)_{\mathcal{U}}$ is weak-weak$^*$ continuous and isometric, and by restricting this representation to the subspace $F$, it is essential so. Indeed, this representation is an essential representation of $L_1(G)$ as following
\begin{align*}
&(\pi^\infty)_{\mathcal{U}}: L_1(G)\rightarrow \mathcal{B}(F),\\
&(\pi^\infty)_{\mathcal{U}}(f)(\xi_n)_{\mathcal{U}}=(\pi(f)\xi_n)_{\mathcal{U}},\quad (\xi_n)_{\mathcal{U}}\in F\subset l_p(\mathbb{N},E)_{\mathcal{U}}.
\end{align*}
Moreover, this representation comes from a representation of $G$ which we still denote it by $((\pi^\infty)_\mathcal{U},F)$.
We need to notice that the space $F$ is a $QSL_p$-space, therefore, it is super-reflexive, and we have 
\begin{align*}
\mathcal{B}(F)=(F^*\widehat{\otimes}F)^* \quad \text{and}\quad F^*\widehat{\otimes}F\subset \mathcal{B}(F)^* .
\end{align*}
Since $(\pi^\infty)_\mathcal{U}:PF_{p,\pi}(G)\rightarrow \mathcal{B}(F)$ is weak-weak$^*$ continuous and isometric, then $(\pi^\infty)_\mathcal{U}^*$ restricted to $F^*\widehat{\otimes}F$ is a quotient map onto $PF_{p,\pi}(G)^*$, so we have
\begin{align*}
PF_{p,\pi}(G)^*\cong F^*\widehat{\otimes}F/\ker(\pi^\infty)_\mathcal{U}^*.
\end{align*}
Indeed, the restricted map $(\pi^\infty)_\mathcal{U}^*:F^*\widehat{\otimes}F\rightarrow PF_{p,\pi}(G)^*$ is of the type of maps in Remark \ref{REMARKPSIMAP}-\eqref{REMARKPSIMAP1}, $\Psi_{p,(\pi^\infty)_\mathcal{U}}$.We note that $\ker(\pi^\infty)_\mathcal{U}^*$ is weak$^*$ closed. Additionally, since $(\pi^\infty)_\mathcal{U}$ is one-to-one, then $\mathcal{R}((\pi^\infty)_\mathcal{U}^*)$, range of $(\pi^\infty)_\mathcal{U}^*$, is weak$^*$ dense in $PF_{p,\pi}(G)^*$, so we have 
$$\overline{A_{p,(\pi^\infty)_\mathcal{U}}{}}^{w^*}=PF_{p,\pi}(G)^*.$$

\end{proof}

\end{proposition}

Next theorem is a generalization of \cite[Theorem 4]{COWLING1979}, in which for a compact subset $K\subset G$, by $A_{p,(\pi^\infty)_{\mathcal{U}}}|_K$ we mean the restriction of functions in $A_{p.(\pi^\infty)_\mathcal{U}}$ to $K$. It is evident that for a function $u\in A_{p,(\pi^\infty)_{\mathcal{U}}}$ and compact set $K\subset G$ we have $\| u|_K\|_{A_{p,(\pi^\infty)_{\mathcal{U}}}|_K}\leq\| u\|_{A_{p,(\pi^\infty)_{\mathcal{U}}}}$, where the norm $\| u|_K\|_{A_{p,(\pi^\infty)_{\mathcal{U}}}|_K}$ is naturally defined to be the infimum of all possible expressions of restricted function $u|_K$ as a coefficient function of representation $((\pi^\infty)_\mathcal{U},F)$.

\begin{theorem}\label{PROPOSITIONDUALOFPFP}
Let $(\pi, E)\in\text{Rep}_p(G)$, then a function $w\in L_\infty(G)$ belongs to $PF_{p,\pi}(G)^*$ with $\|w\|\leq
 C $ if and only if $w|_K\in A_{p,(\pi^\infty)_\mathcal{U}}|_K$ with $\|w|_K\|_{ A_{p,(\pi^\infty)_\mathcal{U}}|_K}\leq C$, for every compact subset $K$ of $G$.
\begin{proof}
\begin{enumerate}
\item[]
 First, we assume that $w|_{K}\in A_{p,(\pi^\infty)_\mathcal{U}}|_K$, with $\| w|_K\|_{A_{p,(\pi^\infty)_\mathcal{U}}|_K}\leq C$, for all compact subset $K\subset G$. Let $f\in L_1(G)$ with compact support $K_f$, then
\begin{align*}
|\langle \pi(f), w\rangle |=|\int_G f(x)w(x) dx|=|\int_{K_f} f(x)w|_{K_f}(x) dx|.
\end{align*}
Since $w|_{K_f}\in A_{p,(\pi^\infty)_\mathcal{U}}|_{K_f}$, then $
w|_{K_f}(\cdot)=\sum_{n=1}^{\infty}\langle (\pi^\infty)_\mathcal{U}(\cdot)\xi_n^f,\eta_n^f\rangle $ with $(\xi_n^f)_{n=1}^\infty\subset F$ and $(\eta_n^f)_{n=1}^\infty\subset F^*$,
so we have
\begin{align*}
|\langle \pi(f), w\rangle |&=|\sum_{n=1}^{\infty}\int_{K_f}f(x)\langle(\pi^\infty)_\mathcal{U}(x)\xi_n^f,\eta_n^f\rangle dx|\\
&=|\sum_{n=1}^{\infty}\int_{K_f}f(x)\langle(\pi^\infty)_\mathcal{U}(x)\xi_n^f,\eta_n^f\rangle dx|\\
&=|\sum_{n=1}^{\infty}\langle(\pi^\infty(f))_\mathcal{U}\xi_n^f,\eta_n^f\rangle |\\
&\leq\|(\pi^\infty(f))_\mathcal{U}\| \sum_{n=1}^{\infty}\|\xi_n^f\|\|\eta_n^f\|\\
&=\|\pi(f)\| \sum_{n=1}^{\infty}\|\xi_n^f\|\|\eta_n^f\|.
\end{align*}
Consequently, we have $|\langle \pi(f), w\rangle |\leq C\| f\|_\pi$, and since compact support functions are dense in $PF_{p,\pi}(G)$, so we have $w\in PF_{p,\pi}(G)^*$, and $\| w\|\leq C$.

\item[]
Now let $ w\in PF_{p,\pi}(G)^*=\overline{A_{p,(\pi^\infty)_\mathcal{U}}{}}^{w^*}$. It follows from Proposition \ref{PROPULTRAAP} and \cite[Lemma 6.5]{RUNDE2005} that there exists a unique $u\in B_p(G)$ which is a coefficient function of the representation $((\pi^\infty)_\mathcal{U},F)$ such that for every $\epsilon >0$ there are vectors $(\xi_n)_{n\in\mathbb{N}}\subset F$ and $(\eta_n)_{n\in\mathbb{N}}\subset F^*$ so that
\begin{align}\label{111}
u(\cdot)=\sum\langle(\pi^\infty)_\mathcal{U}(\cdot)\xi_n,\eta_n\rangle,\quad \|w\|\leq\sum_{n=1}^\infty\|\xi_n\|\|\eta_n\|<\|w\|+\epsilon,
\end{align}
and 
\begin{align}\label{222}
\langle\pi(f),w\rangle=\int_Gf(x)u(x)dx=\langle\pi(f),u\rangle,
\end{align}
which \eqref{111} and \eqref{222} mean that $w$ and $u$ are equal as functionals on $PF_{p,\pi}(G)$, and by Hahn-Banach theorem we have $w=u$ with $\| w\|=\| u\|$. Therefore, the restriction of $w$ to every compact subset $K\subset G$  means the restriction of $u$ to $K$, and we have
\begin{align*}
u|_K\in A_{p,(\pi^\infty)_\mathcal{U}}|_K,\quad \|u|_K\|_{A_{p,(\pi^\infty)_\mathcal{U}}|_K}\leq \| u\|\leq C.
\end{align*}
\end{enumerate}
\end{proof}
\end{theorem}

\begin{remark}\label{REMARKOFDUALOFPFP}
\begin{enumerate}
\item
We follow \cite{ARSAC1976} in notation, and denote $\overline{A_{p,(\pi^\infty)_{\mathcal{U}}}{}}^{w^*}$ by $B_{p,\pi}$, and we call it \textit{$p$-analog of the $\pi$-Fourier-Stieltjes algebra}, which by Proposition \ref{PROPULTRAAP} is the dual space of the space of $p$-pseudofunctions  associated with $(\pi,E)\in\text{Rep}_p(G)$, i.e. the dual space of $PF_{p,\pi}(G)$ through following duality
\begin{align*}
\langle \pi(f),u\rangle=\int_G u(x)f(x)dx,\quad f\in L_1(G),\ u\in B_{p,\pi},
\end{align*}
and as we expect that, we have
\begin{align*}
&\| u\|=\sup_{\|f\|_\pi\leq 1}|\langle \pi(f), u\rangle|=\sup_{\|f\|_\pi\leq 1}|\int_G u(x)f(x)dx|,\quad u\in B_{p,\pi},\\
&\| f\|_\pi=\sup_{\|u\|\leq 1}|\langle \pi(f), u\rangle|=\sup_{\|u\|\leq 1}|\int_G u(x)f(x)dx|,\quad f\in L_1(G).\\
\end{align*}
\item\label{REMARKIDENTITY} It is obvious that $B_{p,\pi}\subset B_p(G)$ is a contractive inclusion for every $(\pi , E)\in\text{Rep}_p(G)$, and if $(\pi,E)$ is a $p$-universal representation it will become an isometric isomorphism.
\item It is valuable to note that if $\mathcal{V}$ is another free ultrafilter as it is described in Proposition \ref{PROPULTRAAP}, then we have
\begin{align*}
\overline{A_{p,(\pi^\infty)_{\mathcal{U}}}{}}^{w^*}=PF_{p,\pi}(G)^*= \overline{A_{p,(\pi^\infty)_{\mathcal{V}}}{}}^{w^*}.
\end{align*}
So, our definition is independent of choosing suitable free ultrafilter, therefore, it is well-defined.
\item\label{REMAPBPMPINC}
For a locally compact group $G$ we have the following relations

\begin{align*}
\overline{A_p(G){}}^{w^*}=B_{p,\lambda_p}\subset B_p(G)\subset \mathcal{M}(A_p(G)),
\end{align*}
and all inclusions will become equalities in the case that $G$ is amenable (See \cite[Theorem 6.6 and Theorem 6.7]{RUNDE2005}).
\end{enumerate}
\end{remark}

\section{Applicable $ p $-completely bounded homomorphisms on $B_p(G)$}\label{APPLICABLEHOMO}

In the following, we study completely boundedness of special type of maps on the $p$-analog of the Fourier-Stieltjes algebras. To provide requirements of forthcoming propositions, Theorem \ref{IMPORTANTMAPPING} plays a critical role. For this aim, we give next lemma that is a kind of application of  Proposition \ref{PROPULTRAAP}.\\
Let $G_0\subset G$, be any subset, and $u:G_0\rightarrow \mathbb{C}$ be a function. By $u^\circ$ we mean
\begin{align*}
u^\circ=\left\{
\begin{array}{ll}
u&\text{on}\; G_0\\
0&\text{o.w.}
\end{array}\right. .
\end{align*}

\begin{lemma}\label{LEMMARESTRICTIONMAP}
Let $G_0$ be an open subgroup of the locally compact group $G$, and $(\pi , E)$ denote the $p$-universal representation of $G$. Then $(\pi_{G_0} ,E)$, restriction of $\pi$ to $G_0$, is a representation of $G_0$, and the restriction mapping of functions in $B_p(G)$ to $G_0$, is a contractive linear homomorphism into $B_{p,\pi_{G_0}}\subset B_p(G_0)$. Moreover, we have the following contractive inclusions
\begin{align*}
 B_{p,\lambda_{p,G_0}}\subset B_{p,\pi_{G_0}}\subset B_p(G_0).
\end{align*}

\begin{proof}

Let us define
\begin{align*}
\pi_{G_0}:G_0\rightarrow \mathcal{B}(E),\quad \pi_{G_0}(x)=\pi(x),\ x\in G_0,
\end{align*}
which obviously  implies that $(\pi_{G_0}, E)\in\text{Rep}_{p}(G_0)$, and by Remark \ref{REMARKGARDELLA} it is contained in a $p$-universal representation of $G_0$, namely $(\rho,F)$, and  we have (up to an isometry)
\begin{align*}
E\subset F,\quad \pi(x)=\pi_{G_0}(x)=\rho(x)|_E,\quad x\in G_0.
\end{align*}
Consequently, through Remark \ref{REMARKOFDUALOFPFP}-(\ref{REMARKIDENTITY}), for a function $u\in B_p(G)$, it can be obtained that
\begin{align*}
u|_{G_0}\in B_{p,\pi_{G_0}}\subset B_p(G_0),\quad \|u|_{G_0}\|_{B_p(G_0)}\leq\| u|_{G_0}\|_{B_{p,\pi_{G_0}}}\leq\|u\|_{B_p(G)}.
\end{align*}
For the inclusion $B_{p,\lambda_{p,G_0}}\subset B_{p,\pi_{G_0}}$, let $u\in A_{p}(G_0)$, then by Remark \ref{REMARKRUNDERUNDE}-\eqref{ITEMFIGARESTRICTION}, we have $u^\circ\in A_p(G)\subset B_p(G)$. Since $B_{p,\pi_{G_0}}$ is the collection of functions in $B_p(G)$ restricted to $G_0$, and since $u=(u^\circ)|_{G_0}$, then we have $u\in B_{p,\pi_{G_0}}$, and it means that $A_{p}(G_0)\subset B_{p,\pi_{G_0}}$. Therefore, via Remark \ref{REMARKOFDUALOFPFP}-\eqref{REMAPBPMPINC}, we have $\overline{A_p(G_0){}}^{w^*}=B_{p,\lambda_{p,G_0}}\subset B_{p,\pi_{G_0}}$. Additionally, one can reach to this inclusion by utilizing Theorem \ref{PROPOSITIONDUALOFPFP}.
\end{proof}
\end{lemma}

\begin{lemma}\label{LEMMAMULTIPLIER}
Let $G_0$ be an open subgroup of the locally compact group $G$, and $u\in\mathcal{M}(A_p(G_0))$. Then we have $u^\circ\in\mathcal{M}(A_p(G))$ with $\|u^\circ\|_{\mathcal{M}(A_p(G))}= \|u\|_{\mathcal{M}(A_p(G_0))}$.
\begin{proof}
Let $u\in \mathcal{M}(A_p(G_0))$ and $v\in A_p(G)$. By the relation $u^\circ\cdot v=(u\cdot v|_{G_0})^\circ$, it can be concluded that $u^\circ \in \mathcal{M}(A_p(G))$, and obviously we have $\|u^\circ\|_{\mathcal{M}(A_p(G))}= \|u\|_{\mathcal{M}(A_p(G_0))}$.
\end{proof}
\end{lemma}

\begin{proposition}\label{PROPEXT}
Let $G$ be a locally compact group and $G_0$ be its open subgroup. Then
\begin{enumerate}
\item
for every $u\in B_p(G_0)$, we have $u^\circ\in\mathcal{M}(A_p(G))$,
\item if $G_0$ is also an amenable subgroup, then for every $u\in B_p(G_0)$, we have $u^\circ\in B_p(G)$.
\end{enumerate}
\begin{proof}
\begin{enumerate}
\item
This part can be concluded by the inclusions in Remark \ref{REMARKOFDUALOFPFP}-\eqref{REMAPBPMPINC} and Lemma \ref{LEMMAMULTIPLIER}.
\item
Since $G_0$ is amenable, then by "Moreover" part in the Lemma \ref{LEMMARESTRICTIONMAP}, and equalities in Remark \ref{REMARKOFDUALOFPFP}-\eqref{REMAPBPMPINC} we have the result.

\end{enumerate}
\end{proof}
\end{proposition}

As an immediate consequence of Proposition \ref{PROPEXT}, we have the next corollary.

\begin{corollary}
Let $G$ and $H$ be locally compact groups, and $\alpha :Y=\cup_{k=1}^nY_k\subset H\rightarrow G$ be a continuous piecewise affine map with disjoint $Y_k\in\Omega_{\text{am-}0}(H)$, for $k=1,\ldots,n$. Then  $u\in B_p(G)$ implies $(u\circ \alpha)^\circ\in B_p(H)$.
\begin{proof}First of all, we note that similar to the case of the Fourier-Stieltjes algebras, a continuous homomorphism $\beta :H\rightarrow G$, between two locally compact groups $G$ and $H$, induces a homomorphism from $B_p(G)$ into $B_p(H)$, by taking $u\in B_p(G)$ to the function $u\circ \beta\in B_p(H)$, see \cite{ILIESPRONK2005}.

Now we divide our proof into two steps. 
\item[Step 1:] First, we let $\alpha: Y=y_0H_0\rightarrow G $ be a continuous affine map, and $\beta :H_0\rightarrow G$ be the homomorphism associated with $\alpha$, as it is explained in Remark \ref{AFFFIINEREMMM}-\eqref{affine-remark}, for an open amenable subgroup $H_0$ of $H$. As we initially explained, the map $u\mapsto u\circ\beta$ is an algebra homomorphism from $B_p(G)$ into $B_p(H_0)$. Consider the following translation maps
\begin{align*}
&L_{{y_0}^{-1}}: B_p(H)\rightarrow B_p(H),\quad L_{{y_0}^{-1}}(u)(h)=u({y_0}^{-1}h),\quad u\in B_p(H),\ h\in H,\\
&L_{\alpha(y_0)}:  B_p(G)\rightarrow B_p(G),\quad L_{y_0}(u)(g)=u(y_0g),\quad u\in B_p(G),\ g\in G,
\end{align*}
then by the following relation, and applying Proposition \ref{PROPEXT}, we have the result
\begin{align*}
(u\circ\alpha)^\circ=L_{{y_0}^{-1}}((L_{\alpha(y_0)}u)\circ\beta)^\circ,\quad u\in B_p(G).
\end{align*}
 
\item[Step 2:] Now let $\alpha: Y\subset H\rightarrow G $ be a continuous piecewise affine map, so by our assumption of amenability, and similar to the Definition \ref{DEFPIECEWISE}, there exist pairwise disjoint sets $Y_k\in\Omega_{\text{am-}0}(H)$, for $k=1,\ldots,n$ with $n\in\mathbb{N}$, and affine maps $\alpha_k:\text{Aff}(Y_k)\subset H\rightarrow G$ such that $Y=\cup_{k=1}^n Y_k$, and $\alpha_k|_{Y_k}=\alpha|_{Y_k}$. By previous step, we know that $(u\circ\alpha_k)^\circ\in B_p(H)$, and since 
\begin{align*}
(u\circ\alpha_k)^\circ=\sum_{k=1}^n(u\circ\alpha_k)^\circ\cdot\chi_{Y_k},
\end{align*}
we have the result via Corollary \ref{COROIDEMCOSETRING}, and the fact that $B_p(H)$ is a Banach algebra.
\end{proof}
\end{corollary}

\begin{remark}\label{REMARKEXTENSION}
\begin{enumerate}
\item\label{REMARKEXTENSIONREP} For an open amenable subgroup $G_0$ of the locally compact group $G$, by Proposition \ref{PROPEXT}, we can say that the space $B_p(G_0)$ is the space of functions which are restriction of functions in $B_p(G)$, those are equal to zero outside of $G_0$. Therefore, while we are working on the $p$-analog of the Fourier-Stieltjes algebras, we may assume that the $p$-universal representation of an open amenable subgroup $G_0$ of $G$ is the restriction of the $p$-universal representation of $G$ to $G_0$.

\item
For an open amenable subgroup $G_0$ of $G$, by Lemma \ref{LEMMARESTRICTIONMAP} and Proposition \ref{PROPEXT}, the restriction mapping from $B_p(G)$ to $B_p(G_0)$ is surjective.
\item
In the case that $G$ is amenable, Proposition \ref{PROPEXT} can be concluded directly from Lemma \ref{LEMMAMULTIPLIER}, via the isometric identification $B_p(G)=\mathcal{M}(A_p(G))$.

\end{enumerate}
\end{remark}

Next theorem is our first main result of this paper, and it will be applied to give the results on weighted homomorphisms on the $p$-analog of the Fourier-Stieltjes algebras. For more clarification, we need to introduce the notion of the $p$-tensor product $E\tilde{\otimes}_p F$ of two $QSL_p$-spaces $E$ and $F$, that is defined in \cite{RUNDE2005}. In fact, Runde introduced the norm $\|\cdot\|_p$ on the algebraic tensor product $E\otimes F$ which benefits from pivotal properties. As an important property of the norm $\|\cdot\|_p$, is the fact that the completion $E\tilde{\otimes}_p F$ of $E\otimes F$ with respect to $\|\cdot\|_p$ is a $QSL_p$-space. Furthermore, for two representations $(\pi,E)$ and $(\rho,F)$ of the locally compact group $G$ in $\text{Rep}_p(G)$, the representation $(\pi\otimes\rho,E\tilde{\otimes}_p F)$ is well-defined and belongs to $\text{Rep}_p(G)$. As a result, for two functions $u(\cdot)=\langle\pi(\cdot)\xi,\eta\rangle$ and $v(\cdot)=\langle\rho(\cdot)\xi',\eta'\rangle$, the pointwise product of them is a coefficient function of the representation $(\pi\otimes\rho, E\tilde{\otimes}_p F)$, i.e. $u\cdot v(\cdot)=\langle(\pi(\cdot)\otimes\rho(\cdot))(\xi\otimes\xi'),\eta\otimes\eta'\rangle$. For more details on $p$-tensor product $\tilde{\otimes}_p$ see \cite[Theorem 3.1 and Corollary 3.2]{RUNDE2005}.

\begin{theorem}\label{IMPORTANTMAPPING} Let $p\in (1,\infty)$ and $G$ be a locally compact group. Then we have the following statements:
\begin{enumerate}
\item
For any $(\pi_p, E_p)\in \text{Rep}_{p}(G)$, the identity map $I:B_{p,\pi_p}\rightarrow B_p(G)$ is a $p$-completely contractive map.
\item
For an open subgroup $G_0$ of $G$, the restriction map $R_{G_0}: B_p(G)\rightarrow B_p(G_0)$, is a $p$-completely contractive homomorphism.
\item
For an element $a\in G$, the translation map $L_a:B_p(G)\rightarrow B_p(G)$, defined through $L_a(u)={}_a u$, where ${}_au(x)=u(ax)$, for $x\in G$, is a $p$-completely contractive map.
\item
For a closed normal subgroup $G_1$ of $G$, let $q:G\rightarrow G/G_1$ be the canonical quotient map. Then the homomorphism $\Phi_q: B_p(G/G_1)\rightarrow B_p(G)$, with $\Phi_q(u)= u\circ q$, is a $p$-completely contractive homomorphism.
\item
For an open amenable subgroup $G_2$ of $G$, the extension map $E_{G_2}:B_p(G_2)\rightarrow B_p(G)$ is a $p$-completely contractive homomorphism.
\item For an open coset $Y$ of an open subgroup $G_2$ of $G$, the map $M_Y:B_p(G)\rightarrow B_p(G)$, with $M_Y(u)=u\cdot\chi_Y$, is $p$-completely contractive homomorphism. More generally, for a set $Y\in\Omega_0(G)$, the map $M_Y$ is a $p$-completely bounded homomorphism.
\end{enumerate}

\begin{proof}

\begin{enumerate}

\item\label{ITEMIDENTITYMAP}
We want to prove that for each $(\pi_p,E_p)\in\text{Rep}_{p}(G)$, the following map is a $p$-complete contraction.
\begin{align}\label{ID}
I: B_{p,\pi_p}\rightarrow B_p(G),\quad I(u)=u .
\end{align}

Let $(\pi ,E)$ be a $p$-universal representation of $G$ that contains the representation $(\pi_p,E_p)$. Following relations hold between $(\pi_p , E_p)$, and $(\pi, E)$.
\begin{align*}
E_p\subset E, \quad \pi_p(x)=\pi(x)|_{E_p},\quad\text{and}\quad\pi_p(f)=\pi(f)|_{E_p},\quad x\in G,\; f\in L_1(G).
\end{align*}
Since $\pi_p(f)=\pi(f)|_{E_p}$, then $\|\pi_p(f)\|\leq\|\pi(f)\|$. Additionally, the map $I$ is weak$^*$-weak$^*$ continuous, and it is a contraction by \cite[Theorem 6.6-$(i)$]{RUNDE2005}.
Define
\begin{align*}
{}_*I: UPF_{p}(G)\rightarrow PF_{p,\pi_p}(G),\quad {}_*I(\pi(f))=\pi(f)|_{E_p}=\pi_p(f),
\end{align*}
then ${}_*I$ is the predual of the map \eqref{ID}. Because, we have $\langle\pi(f), I(u)\rangle=\langle\pi_p(f),u\rangle$, for every $f\in L_1(G)$ and $u\in B_{p,\pi_p}$. Following calculations indicate that ${}_*I$ is a $p$-complete contraction: for each $n\in\mathbb{N}$, and $(\pi_p(f_{ij}))\in\mathbb{M}_n(PF_{p,\pi_p}(G))$ we have
\begin{align*}
\|(\pi_{p,\pi_p}(f_{ij}))\|_n&=\sup\{\|(\pi_p(f_{ij}))(\xi_{j})_{j=1}^n\|\ :\  (\xi_j)_{j=1}^n\in\mathbb{M}_n(E_p),\; \sum_{j=1}^n\|\xi_j\|^p\leq 1\}\\
&=\sup\{\|(\pi(f_{ij}))(\xi_{j})_{j=1}^n\|\ :\  (\xi_j)_{j=1}^n\in\mathbb{M}_n(E_p),\; \sum_{j=1}^n\|\xi_j\|^p\leq 1\}\\
&\leq\sup\{\|(\pi(f_{ij}))(\xi_{j})_{j=1}^n\|\ :\  (\xi_j)_{j=1}^n\in\mathbb{M}_n(E),\; \sum_{j=1}^n\|\xi_j\|^p\leq 1\}\\
&=\|(\pi(f_{ij}))\|_n,
\end{align*}
so we have $\|(\pi_p(f_{ij}))\|_n\leq \|(\pi(f_{ij}))\|_n$, and by this, it is concluded that
\begin{align*}
\|I\|_{\text{p-cb}}\leq \|{}_*I\|_{\text{p-cb}}\leq 1 .
\end{align*} 

\item\label{ITEMRESTRICTIONMAP}
Let $G_0\subset G$ be an open subgroup and consider the following map:
\begin{align*}
R_{G_0}: B_p(G)\rightarrow B_p(G_0),\quad R_{G_0}(u)=u|_{G_0}.
\end{align*}
Let $(\pi, E)$ be a $p$-universal representation of $G$, and $(\pi_{G_0},E)$ be the restriction of $(\pi,E)$ to $G_0$. Their liftings are related as following
\begin{align}\label{RESREQ}
\pi_{G_0}(f)=\pi(f^\circ),\quad f\in L_1(G_0).
\end{align}
In addition, range of the map $R_{G_0}$ is the space $B_{p,\pi_{G_0}}\subset B_p(G_0)$, as it is described in Lemma \ref{LEMMARESTRICTIONMAP}. This map is weak$^*$-weak$^*$ continuous by the relation \eqref{RESREQ}, and 
\begin{align*}
\langle \pi_{G_0}(f),u|_{G_0}\rangle=\langle \pi(f^\circ),u\rangle,\quad f\in L_1(G_0),\; u\in B_p(G).
\end{align*}
So, we may define ${}_*R_{G_0}$ as following
\begin{align*}
{}_*R_{G_0}:PF_{p,\pi_{G_0}}(G_0)\rightarrow UPF_{p}(G),\quad {}_*R_{G_0}(\pi_{G_0}(f))=\pi(f^\circ),\quad f\in L_1(G_0).
\end{align*}

We have
\begin{align*}
\langle\pi_{G_0}(f), ({}_*R_{G_0})^*(u)\rangle&=\langle {}_*R_{G_0}(\pi_{G_0}(f)),u\rangle\\
&=\langle \pi(f^\circ),u\rangle\\
&=\int_G u(x) f^\circ(x)dx\\
&=\int_G u|_{G_0}(x) f(x)dx\\
&=\langle \pi_{G_0}(f),u|_{G_0}\rangle\\
&=\langle \pi_{G_0}(f),R_{G_0}(u)\rangle.
\end{align*}
Therefore, $({}_*R_{G_0})^*=R_{G_0}$. Additionally, by \eqref{RESREQ} we have ${}_*R_{G_0}(\pi_{G_0}(f))=\pi(f^\circ)=\pi_{G_0}(f)$, so ${}_*R_{G_0}$ is an identity map which is $p$-completely isometric
\begin{align*}
\|{}_*R^{(n)}(\pi_{G_0}(f_{ij}))\|_n=
\|(\pi(f_{ij}^\circ))\|_n=\|(\pi_{G_0}(f_{ij}))\|_n,
\end{align*}
therefore, $\|R_{G_0}\|_{\text{p-cb}}=\|({}_*R_{G_0})^*\|_{\text{p-cb}}
\leq \|{}_*R_{G_0}\|_{\text{p-cb}}=1$.

\item\label{ITEMTRANSLATIONMAP}
Now we want to prove that for $a\in G$, the following map is a $p$-complete isometry
\begin{align*}
L_a:B_p(G)\rightarrow B_p(G),\quad L_a(u)={}_au,\quad {}_au(x)=u(ax),\; x\in G.
\end{align*}
Predual of the map $L_a$ is as following
\begin{align*}
&{}_*L_a :UPF_{p}(G)\rightarrow UPF_p(G),\quad {}_*L_a(\pi(f))=\pi(\lambda_p(a)f)),
\end{align*}
and it is clearly $p$-completely contractive, and consequently, this is true for $L_a$. On the other hand, the map $L_a$ has the inverse $L_{a^{-1}}$, and similar to $L_a$, it is $p$-completely contractive which makes $L_a$ to be $p$-completely isometric map.

\item\label{ITEMQUOTIENTMAP} Let $G_1\subset G$ be a closed normal subgroup. Let 
\begin{align*}
q:G\rightarrow G/G_1,\quad q(x)=xG_1,\quad x\in G,
\end{align*}
be the canonical quotient map, and 
\begin{align*}
\Phi_q :B_p(G/G_1)\rightarrow B_p(G),\quad\Phi_q(u)=u\circ q .
\end{align*}
Let $(\rho, F)$ be a $p$-universal representation of $G/G_1$. Then obviously we have $(\rho\circ q, F)\in\text{Rep}_{p}(G)$, and this representation is contained in a $p$-universal representation $(\pi, E)$ of $G$ which implies that (up to an isometry)
\begin{align*}
& F\subset E,\quad \rho\circ q(x)=\pi(x)|_{F},\quad\rho\circ q(f)=\pi(f)|_{F},\quad x\in G,\ f\in L_1(G).
\end{align*}

Let us define closed subspace $K$ of $E$, which itself is a $QSL_p$-space, through
\begin{align*}
K=\{\xi\in E\ :\ \pi(x)\xi=\xi,\; \forall\; x\in G_1\},
\end{align*}
and consider the representation of $G$, for which an element $x\in G$ goes to the restriction of $\pi(x)$ to $K$, so we can induce the following representation for $G/G_1$
\begin{align*}
\tilde{\pi}:G/G_1\rightarrow\mathcal{B}(K),\quad\tilde{\pi}(xG_1)=\pi(x)|_K,
\end{align*}
by the definition of $K$, we have $F\subset K$, therefore, $(\rho , F)\subset (\tilde{\pi}, K)$. It is obtained that every $p$-universal representation of $G/G_1$, like $(\rho, F)$, is contained in a representation of $G/G_1$, as we described, $(\tilde{\pi},K)$, that is induced by the $p$-universal representation $(\pi, E)$ of $G$. So, we can work by $(\tilde{\pi},K)$ as the $p$-universal representation of $G/G_1$.

Now, let us consider the following map
\begin{align*}
\Phi_q: B_{p}(G/G_1)\rightarrow B_{p,\tilde{\pi}\circ q}\subset B_p(G),\quad\Phi_q(u)=u\circ q,
\end{align*}
which is at least a contractive isomorphism into the subalgebra of $B_p(G)$, of functions which are constant on each coset of $G_1$. For functions $f\in L_1(G)$, and $u\in B_p(G/G_1)$, we have
\begin{align}\label{EQPHQ}
\langle\pi(f),\Phi_q(u)\rangle=\langle\pi(f),u\circ q\rangle=\langle\tilde{\pi}(Pf),u\rangle,
\end{align}
where the map $P:L_1(G)\rightarrow L_1(G/G_1)$ is defined \cite{FOLLAND1995}:
\begin{align*}
Pf(xG_1)=\int_{G_1} f(xg)dg,\qquad f\in C_c(G).
\end{align*}
This implies that the map $\Phi_q$ is weak$^*$-weak$^*$ continuous, and by this we define the predual map ${}_*\Phi_q$, as following:
\begin{align*}
{}_*\Phi_q : PF_{p,\tilde{\pi}\circ q}(G)\rightarrow UPF_{p}(G/G_1),\quad {}_*\Phi_q(\tilde{\pi}\circ q(f))=\tilde{\pi}(P f),\; f\in L_1(G),
\end{align*}
which by \eqref{EQPHQ} we have $({}_*\Phi_q)^*=\Phi_q$. For a function $v:G\rightarrow\mathbb{C}$ that is constant on the cosets of $G_1$, by $\tilde{v}$ we denote
\begin{align*}
\tilde{v}:G/G_1\rightarrow\mathbb{C},\quad \tilde{v}(xG_1)=v(x),\quad x\in G.
\end{align*}
We need to note that for $f\in L_1(G)$, $\xi\in K$ and $\eta\in K^*$, we have
\begin{align*}
\langle\tilde{\pi}\circ q(f)\xi,\eta\rangle&=\int_{G}f(x)\underbrace{\langle\tilde{\pi}\circ q(x)\xi,\eta\rangle}_{v(x)}dx\\
&=\int_{G}f(x)v(x)dx\\
&=\int_{G/G_1}P(f\cdot v)(xG_1)dxG_1\\
&=\int_{G/G_1}\tilde{v}(xG_1)Pf(xG_1)dxG_1\\
&=\int_{G/G_1}Pf(xG_1)\langle\tilde{\pi}(xG_1)\xi,\eta\rangle dxG_1\\
&=\langle\tilde{\pi}(Pf)\xi,\eta\rangle,
\end{align*}
so we have $\tilde{\pi}\circ q(f)=\tilde{\pi}(Pf)$, which means that the predual map ${}_*\Phi_q$ is an identity map that is $p$-completely isometric map via the following computation
\begin{align*}
\|{}_*\Phi_q^{(n)}(\tilde{\pi}\circ q(f_{ij}))\|_n=
\|(\tilde{\pi}(Pf_{ij}))\|_n=
\|(\tilde{\pi}\circ q(f_{ij}))\|_n.
\end{align*}
Therefore, we have $\|\Phi_q\|_{\text{p-cb}}\leq 1$.

\item\label{ITEMEXTENSIONMAP} 
Let $G_2\subset G$, be an open amenable subgroup, and $u\in B_p(G_2)$. Since by Proposition \ref{PROPEXT} we have $u^\circ \in B_p(G)$, then we are allowed to define
\begin{align*}
&E_{G_2}: B_p(G_2)\rightarrow B_p(G),\quad E_{G_2}(u)=u^\circ .
\end{align*}
Let $(\pi ,E)$ be a $p$-universal representation of $G$. We denote the restriction of $(\pi,E)$ to $G_2$ by $(\pi_{G_2}, E)$ which is a $p$-universal representation of $G_2$ via Remark \ref{REMARKEXTENSION}-\eqref{REMARKEXTENSIONREP}. 
We note that by the relation
\begin{align}\label{EQEXT}
\langle \pi(f), u^\circ\rangle=\langle\pi_{G_2}(f|_{G_2}),u\rangle,\quad f\in L_1(G),\ u\in B_p(G_2),
\end{align}
the map $E_{G_2}$ is weak$^*$-weak$^*$ continuous, so we define the predual map ${}_*E_{G_2}$, as following:
\begin{align*}
{}_*E_{G_2}: UPF_p(G)\rightarrow UPF_p(G_2),\quad {}_*E_{G_2}(\pi(f)):=\pi_{G_2}(f|_{G_2}),
\end{align*}
which by \eqref{EQEXT} we have $({}_*E_{G_2})^*=E_{G_2}$.
We need to take notice of the fact that since $\chi_{G_2}\in B_p(G)$, via Theorem \ref{IDMTHEOREM}-\eqref{idm3}, $\chi_{G_2}$ is a normalized coefficient function of $(\pi, E)$, i.e. there are $\xi_\chi\in E$, and $\eta_\chi\in E^*$ with $\|\xi_\chi\|=\|\eta_\chi\|=1$ so that $\chi_{G_2}(\cdot)=\langle\pi(\cdot)\xi_\chi,\eta_\chi\rangle$. Also, for $g\in L_1(G_2)$, and $\xi\in E$, and $\eta\in E^*$, we have
\begin{align*}
\langle\pi_{G_2}(g)\xi,\eta\rangle=\langle\pi(g^\circ)\xi,\eta\rangle,
\end{align*}
and for $f\in L_1(G)$, $\xi\in E$, and $\eta\in E^*$ we have
\begin{align}\label{EQ1}
\langle\pi_{G_2}(f|_{G_2})\xi,\eta\rangle=\langle\pi(f\chi_{G_2})\xi,\eta\rangle.
\end{align}
On the other hand,
\begin{align*}
\langle\pi(f\chi_{G_2})\xi,\eta\rangle&=\int_{G}f(x)\chi_{G_2}(x)\langle\pi(x)\xi,\eta\rangle dx\\
&=\int_{G}f(x)\langle\pi(x)\xi_\chi,\eta_\chi\rangle\langle\pi(x)\xi,\eta\rangle dx\\
&=\int_{G}f(x)\langle(\pi(x)\otimes\pi(x))(\xi_\chi\otimes\xi),\eta_\chi\otimes\eta\rangle dx\\
&=\langle(\pi\otimes\pi(f))(\xi_\chi\otimes\xi),\eta_\chi\otimes\eta\rangle.
\end{align*}
Therefore, by combining last equality with \eqref{EQ1}, we have
\begin{align}\label{EQNEW}
\langle\pi_{G_2}(f|_{G_2})\xi,\eta\rangle=\langle(\pi\otimes\pi(f))(\xi_\chi\otimes\xi),\eta_\chi\otimes\eta\rangle,\quad f\in L_1(G),\ \xi\in E,\ \eta\in E^*.
\end{align}
Additionally, since $(\pi,E)$ is a $p$-universal representation, and we have
\begin{align*}
(\pi, E)\subset (\pi\otimes\pi, E\tilde{\otimes}_p E),
\end{align*}
thus $(\pi\otimes\pi, E\tilde{\otimes}_p E)$ can be assumed as a $p$-universal of $G$. Let
\begin{align*}
{}_*E_{G_2}^{(n)}: \mathbb{M}_n(UPF_p(G))\rightarrow \mathbb{M}_n(UPF_p(G_2)),\quad {}_*E_{G_2}^{(n)}(\pi(f_{ij})):=(\pi_{G_2}(f_{ij}|_{G_2})),
\end{align*}
then via \eqref{EQNEW} we have
\begin{align*}
\|{}_*E_{G_2}^{(n)}(\pi(f_{ij}))\|_n^p&=\|(\pi_{G_2}(f_{ij}|_{G_2}))\|_n^p\\
&=\sup\{|\sum_{i,j=1}^n\langle\pi_{G_2}(f_{ij}|_{G_2})\xi_j,\eta_i\rangle|\ :\ \sum_{j=1}^n\|\xi_j\|^p\leq 1,\;  \sum_{i=1}^n\|\eta_i\|^{p'}\leq 1\}\\
&=\sup\{|\sum_{i,j=1}^n\langle(\pi\otimes\pi(f_{ij}))(\xi_j\otimes\xi_\chi),(\eta_i\otimes\eta_\chi\rangle|\ :\ \sum_{j=1}^n\|\xi_j\|^p\leq 1,\;  \sum_{i=1}^n\|\eta_i\|^{p'}\leq 1\}\\
&\leq\|(\pi\otimes\pi(f_{ij}))\|_n^p,
\end{align*}
and since norm of $UPF_p(G)$ is independent of choosing $p$-universal representation then we have $\|{}_*E_{G_2}\|_{\text{p-cb}}\leq 1$, which implies that $\|E_{G_2}\|_{\text{p-cb}}\leq 1$.

\item\label{ITEMMULTIPLICATIONYYY}
By Corollary \ref{COROIDEMCOSETRING}, the map $M_Y: B_p(G)\rightarrow B_p(G)$ with $M_Y(u)=u\cdot \chi_Y$ is well-defined, and 
\begin{align*}
\|M_Y\|\leq 2^{m_Y} .
\end{align*}
On the other hand, by the following relation this map is weak$^*$-weak$^*$ continuous
\begin{align}\label{EQMULTBP}
\langle\pi(f),u\cdot\chi_Y\rangle=\langle\pi(f\cdot\chi_Y),u\rangle,\quad f\in L_1(G),\ u\in B_p(G).
\end{align}
So, one may define its predual map as following
\begin{align*}
{}_*M_Y:UPF_p(G)\rightarrow UPF_p(G),\quad {}_*M_Y(\pi(f))=\pi(f\cdot\chi_Y),
\end{align*} 
and by \eqref{EQMULTBP} we have $({}_*M_Y)^*=M_Y$.\\
Step 1: For proving the claim, first we let $Y$ be an open coset itself. By Theorem \ref{IDMTHEOREM}-\eqref{idm3}, the function $\chi_Y$ is a normalized coefficient function of representation $(\pi,E)$ which means  that there are elements $\xi_Y\in E$, and $\eta_Y\in E^*$ with $\|\xi_Y\|=\|\eta_Y\|=1$ such that
\begin{align*}
\chi_Y(\cdot)=\langle\pi(\cdot)\xi_Y,\eta_Y\rangle .
\end{align*}
So, for a matrix $(\pi(f_{ij}))\in \mathbb{M}_n(UPF_p(G))$, we have
{\small \begin{align*}
\|(\pi(f_{ij}\cdot\chi_Y))\|_n&=\sup\{|\sum_{i,j=1}^n\langle\pi(f_{ij}\cdot\chi_Y)\xi_j,\eta_i\rangle|\ :\ \sum_{j=1}^n\|\xi_j\|^p\leq 1,\; \sum_{i=1}^n\|\eta_i\|^p\leq 1 \}\\
&=\sup\{|\sum_{i,j=1}^n\int_G f_{ij}(x)\chi_Y(x)\langle\pi(x)\xi_j,\eta_i\rangle dx|\ :\ \sum_{j=1}^n\|\xi_j\|^p\leq 1,\; \sum_{i=1}^n\|\eta_i\|^p\leq 1 \}\\
&=\sup\{|\sum_{i,j=1}^n\int_G f_{ij}(x)\langle(\pi(x)\otimes\pi(x))(\xi_j\otimes\xi_Y),\eta_i\otimes\eta_Y\rangle dx|\ :\ \sum_{j=1}^n\|\xi_j\|^p\leq 1,\; \sum_{i=1}^n\|\eta_i\|^p\leq 1 \}\\
&\leq \sup\{|\sum_{i,j=1}^n\langle\pi\otimes\pi(f_{ij})\phi_j,\psi_i\rangle|\ :\ \sum_{j=1}^n\|\phi_j\|_{E\tilde{\otimes}_pE}^p\leq 1,\; \sum_{i=1}^n\|\psi_i\|_{E^*\tilde{\otimes}_{p'}E^*}^p\leq 1 \}\\
&=\|(\pi\otimes\pi(f_{ij}))\|_n.
\end{align*}}
By these computations, we obtain that the map ${}_*M_Y$ is a $p$-complete contraction. Therefore, we have $\|M_Y\|_{\text{p-cb}}\leq 1$. Note that in the above calculations, we used a relation similar to \eqref{EQNEW} and  an argument about independence of choosing $p$-universal representation.\\
Step 2: Now for $Y=Y_0\backslash\cup_{i=1}^m Y_i \in\Omega_0(G)$, from \eqref{EQCHI} we have,
\begin{align*}
M_{Y}=M_{{Y_0}}-(\sum_{i=1}^m M_{Y_i} -\sum_{i,j}M_{{Y_i\cap Y_j}}+\sum_{i,j,k}M_{{Y_i\cap Y_j\cap Y_k}}+\ldots+(-1)^{m+1}M_{{Y_1\cap\ldots Y_m}}).
\end{align*}
Therefore, we have $\|M_Y\|_{\text{p-cb}}\leq 2^{m_Y}$.

\end{enumerate}
\end{proof}
\end{theorem}

\begin{remark}
\begin{enumerate}
\item
The importance of Theorem \ref{IMPORTANTMAPPING}-\eqref{ITEMIDENTITYMAP} is that while we are working with maps with ranges as subspaces of the $p$-analog of the Fourier-Stieltjes algebras, we just need to restrict ourselves to their ranges, as what we have done in the rest of Theorem \ref{IMPORTANTMAPPING}. 
\item
In the proof of Theorem \ref{IMPORTANTMAPPING}-\eqref{ITEMMULTIPLICATIONYYY}, if $Y=y_2G_2$, for an open amenable subgroup $G_2$ of $G$, and some $y_2\in G$,  then  by Theorem \ref{IMPORTANTMAPPING}-\eqref{ITEMRESTRICTIONMAP}-\eqref{ITEMTRANSLATIONMAP}-\eqref{ITEMEXTENSIONMAP} we can find out that the map $M_Y$ is a $p$-complete contraction through the following relation
\begin{align*}
M_Y=L_{{y_2}^{-1}}\circ E_{G_2}\circ R_{G_2}\circ L_{y_2}.
\end{align*}

\end{enumerate}
\end{remark}

\section{$p$-Completely homomorphisms on $B_p(G)$ induced by proper  piecewise affine maps}\label{WEIGHTEDHOMO}

As an application of previous sections, we are ready to study on homomorphisms $\Phi_\alpha : B_p(G)\rightarrow B_p(H)$ of the form
\begin{align*}
\Phi_\alpha u= \left\{
\begin{array}{ll}
u\circ \alpha & \text{on} \: Y\\
0 & \text{o.w}
\end{array}\right.,\quad u\in B_p(G),
\end{align*}
for the proper continuous piecewise affine map $\alpha :Y\subset H\rightarrow G$ with $Y=\cup_{i=1}^nY_i$ and $Y_i\in\Omega_{\text{am-}0}(H)$, which are pairwise disjoint, for $i=1,\ldots,n$. We will give some results in the sequel. For our aim we need the following lemma. For general form of this lemma, see \cite[Lemma 1]{ILIE2013}, and related references there.
\begin{lemma}\label{LEMMACALPHA}
Let $G$ and $H$ be locally compact groups and $\alpha :H\rightarrow G$ be a proper homomorphism that is onto, then there is a constant $c_\alpha>0$, such that
\begin{align*}
\int_H f\circ\alpha(h)dh=c_\alpha\int_G f(x)dx,\quad f\in L_1(G).
\end{align*}

\end{lemma}

\begin{proposition}\label{ITEMHOMOMORPHISMPHI}
Let $G$ and $H$ be locally compact groups and $\alpha :H\rightarrow G$ be a proper continuous group homomorphism. Then the homomorphism $\Phi_\alpha : B_p(G)\rightarrow B_p(H)$, of the form $\Phi_\alpha(u)=u\circ\alpha$, is well-defined and $p$-completely contractive homomorphism.

\begin{proof}

Let $(\pi,E)$ be a $p$-universal representation of $G$. Obviously, $(\pi\circ\alpha , E)\in\text{Rep}_p(H)$, and $\Phi_\alpha$ is a contractive map so that its range is the subspace of $B_p(H)$ of functions which are coefficient functions of the representation $(\pi\circ\alpha , E)$. We will divide our proof into two steps.

\begin{enumerate}
\item[Step 1:]
First, we suppose that $\alpha :H\rightarrow G$ is a continuous isomorphism.
In this case, $(\pi\circ\alpha , E)$ is a $p$-universal representation of $H$, and by Lemma \ref{LEMMACALPHA},  for every $f\in L_1(H)$ and $u\in B_p(G)$, we have
\begin{align*}
\langle\pi\circ\alpha(f),u\circ\alpha\rangle&=\int_H f(h)u\circ\alpha(h)dh\\
&=\int_H (f\circ\alpha^{-1})\circ\alpha(h)u\circ\alpha(h)dh\\
&=c_\alpha\int_G f\circ\alpha^{-1}(x)u(x)dx\\
&=c_\alpha\langle\pi(f\circ\alpha^{-1}),u\rangle.
\end{align*}
Consequently, the map $\Phi_\alpha$ is weak$^*$-weak$^*$ continuous, and we define
\begin{align*}
{}_*\Phi_\alpha :UPF_p(H)\rightarrow UPF_p(G),\quad{}_*\Phi_\alpha(\pi\circ\alpha(f)):=c_\alpha\pi(f\circ\alpha^{-1}).
\end{align*}
According to the above relation, we have $({}_*\Phi_\alpha)^*=\Phi_\alpha$. On the other hand, for every $\xi\in E$ and $\eta\in E^*$, we have
\begin{align*}
\langle\pi\circ\alpha(f)\xi,\eta\rangle & =\int_{H}f(h)\langle\pi\circ\alpha(h)\xi,\eta\rangle dh\\
& =\int_{H}f\circ\alpha^{-1}\circ\alpha(h)\langle\pi\circ\alpha(h)\xi,\eta\rangle dh\\
& =c_\alpha\int_{G}f\circ\alpha^{-1}(x)\langle\pi(x)\xi,\eta\rangle dx\\
&=\langle c_\alpha\pi(f\circ\alpha^{-1})\xi,\eta\rangle,
\end{align*}
which means $\pi\circ\alpha (f)=c_\alpha\pi(f\circ\alpha^{-1})$.
Consequently, ${}_*\Phi_\alpha$ is an identity map, so is a $p$-complete isometry
\begin{align*}
\|{}_*\Phi_\alpha^{(n)}(\pi\circ\alpha(f_{i,j}))\|_n & =\|(c_\alpha\pi(f_{i,j}\circ\alpha^{-1}))\|_n =\|(\pi\circ\alpha(f_{i,j}))\|_n.
\end{align*}
Therefore, $\|\Phi_\alpha\|_{\text{p-cb}}\leq\|{}_*\Phi_\alpha\|_{\text{p-cb}}=1$.

\item[Step 2:]
Now let $\alpha :H\rightarrow G$ be any proper continuous homomorphism. Let $G_0=\alpha(H)$, and $N=\ker\alpha$. Let us define
\begin{align*}
\tilde{\alpha} :H/N\rightarrow G_0,\quad \tilde{\alpha}(xN)=\alpha(x),
\end{align*}
then by Proposition \ref{ramirez}, the map $\tilde{\alpha}$ is a continuous isomorphism, $N$ is a compact normal subgroup of $H$, and $G_0$ is an open subgroup of $G$. Therefore, $\alpha=\tilde{\alpha}\circ q $. By Step 1,  the map $\Phi_{\tilde{\alpha}}$ is $p$-completely contractive, and because of the following composition, $\Phi_\alpha$ is $p$-completely contractive, via Theorem \ref{IMPORTANTMAPPING}-\eqref{ITEMRESTRICTIONMAP}-\eqref{ITEMQUOTIENTMAP}.
\begin{align*}
\Phi_\alpha=\Phi_q\circ\Phi_{\tilde{\alpha}}\circ R_{G_0} .
\end{align*}

\end{enumerate}

\end{proof}

\end{proposition}
For the next proposition, we have to put the amenability assumption on the subgroups of $H$, because of Proposition \ref{PROPEXT}.

\begin{proposition}\label{ITEMAFFINE}
Let $G$ and $H$ be two locally compact groups, $Y$ be an open coset of an open amenable subgroup of $H$, and $\alpha :Y\subset H\rightarrow G$ be a continuous proper affine map. Then the map $\Phi_\alpha : B_p(G)\rightarrow B_p(H)$, defined as
\begin{align*}
\Phi_\alpha(u)=\left\{
\begin{array}{ll}
u\circ\alpha, & \text{on}\; Y,\\
0,&\text{o.w.}
\end{array}\right.,\quad u\in B_p(G),
\end{align*}
is $p$-completely contractive. More generally, if $\alpha$ is a continuous proper piecewise affine map, and $Y=\cup_{i=1}^nY_i$, where disjoint sets $Y_i$ belong to $\Omega_{\text{am-}0}(H)$, then the map $\Phi_\alpha$ is $p$-completely bounded.

\begin{proof}

Let $\alpha :Y=y_0H_0\rightarrow G$ be a continuous proper affine map on the open coset $Y=y_0H_0$, and $H_0$ be an open amenable subgroup of $H$, for which by Remark \ref{AFFFIINEREMMM}-\eqref{affine-remark}, there exists a continuous group homomorphism $\beta :H_0\subset H\rightarrow G$ associated to $\alpha$ such that
\begin{align*}
\beta(h)=\alpha(y_0)^{-1}\alpha(y_0h),\quad h\in H_0.
\end{align*}
which is proper via Remark \ref{AFFFIINEREMMM}-\eqref{HOMAFFPROPER}. Now consider the following composition
\begin{align*}
\Phi_\alpha= L_{{y_0}^{-1}}\circ E_{H_0}\circ\Phi_\beta\circ L_{\alpha(y_0)},
\end{align*}
 then by Proposition \ref{ITEMHOMOMORPHISMPHI}, and Theorem \ref{IMPORTANTMAPPING}-\eqref{ITEMTRANSLATIONMAP}-\eqref{ITEMEXTENSIONMAP} the map $\Phi_\alpha$ is $p$-completely contractive homomorphism.

Next, we consider the piecewise affine case. Let the map $\alpha :Y\subset H\rightarrow G$ be a continuous piecewise affine map. Then for some $n\in\mathbb{N}$, and $i=1,\ldots,n$, there are disjoint sets $Y_i\in\Omega_{\text{am-}0}(H)$, such that $Y=\cup_{i=1}^{n}Y_i$, and $\alpha_i:{Aff(Y_i)}\rightarrow G$ which are affine maps, and $\alpha_i|_{Y_i}=\alpha|_{Y_i}$. Additionally, by Remark \ref{AFFFIINEREMMM}-\eqref{LEMMA8888}, each affine map $\alpha_i$ is proper. Therefore, by considering
\begin{align*}
\Phi_\alpha= \sum_{i=1}^n M_{Y_i}\circ\Phi_{\alpha_i},
\end{align*}
and through the above computations for the maps $\Phi_{\alpha_i}$, we have
\begin{align*}
\|\Phi_\alpha\|_{\text{p-cb}}\leq \sum_{i=1}^n 2^{m_{Y_i}},
\end{align*}
where $m_{Y_i}$ is the corresponding number to each $Y_i$, as it is in Theorem \ref{IMPORTANTMAPPING}-\eqref{ITEMMULTIPLICATIONYYY}.
\end{proof}
\end{proposition}


\begin{thebibliography}{99}


\bibitem{ARSAC1976}
G. Arsac,
{\it Sur l'espace de Banach engendr\'e par les coefficients d'une repr\'esentation unitaire},
  Publ. D\'ep. Math. (Lyon) \textbf{13} (1976) 1-101.



\bibitem{BLECHER1992}
D.P. Blecher,
{\it The standard dual of an operator space},
Pacific Math. J. \textbf{153} (1992) 15-30.






\bibitem{COHEN1960-1}
 P.J. Cohen,
{\it On a conjecture of Littlewood and idempotent measures},
Amer. J. Math. \textbf{82} (1960) 191-212.


\bibitem{COHEN1960-2}
 P.J. Cohen,
{\it On homomorphisms of group algebras},
Amer. J. Math. \textbf{82} (1960) 213-226.



\bibitem{COWLING1979}
 M. Cowling,
{\it An application of Littlewood-Paley theory in harmonic analysis},
Math. Ann. \textbf{241}(1) (1979) 83-96.
















\bibitem{DAWS2010}
M. Daws,
{\it $p$-Operator spaces and Fig\`a-Talamanca-Herz algebras},
J. Operator Theory \textbf{63} (1) (2010) 47-83.








\bibitem{DUNKLRAMIREZ1971}
C.F. Dunkl, D.E. Ramirez,
{\it Homomorphisms on groups and induced maps on certain algebras of measures},
Trans. Amer. Math. Soc. \textbf{160} (1971) 475-485. 














\bibitem{EFFROSRUAN1991}
E.G. Effros, Z.-J. Ruan,
{A new approach to operator spaces},
Canad. J. Math. \textbf{34} (1991) 329-337.



\bibitem{EYMARD1964}
 P. Eymard,
{\it L'alg\'ebre de Fourier d'un groupe localement compact},
Bull. Soc. Math. France \textbf{92} (1964) 181-236.



\bibitem{FHHMPZ2001}
 M. Fabian, P. Habala, P. H\'ajek, V. Montesinos Santalucia, J. Pelant, V. Zizler,
 {\it Functional Analysis and Infinite Dimensional Geometry},
 CMS Books Math. \textbf{8} Springer-Verlag, New York, 2001.



\bibitem{FIGATALAMANCA1965}
 A. Fig\`a-Talamanca,
{\it Translation invariant operators in $L_p$},
Duke Math. J. \textbf{32} (1965) 495-501.






\bibitem{FOLLAND1995}
G. B. Folland,
{\it A Course in Abstract Harmonic Analysis},
CRC Press, Boca Raton, Fla. 1995.






\bibitem{FORREST1994}
B.E. Forrest, 
{\it Amenability and the structure of the algebras $A_p(G)$}, 
Trans. Amer. Math. Soc. \textbf{343} (1994) 233-243.
 
\bibitem{GARDELLATHIEL2015}
E. Gardella, H. Thiel,
{\it Group algebras acting on $L_p$-spaces},  J. Fourier Anal. Appl. \textbf{21} (6) (2015) 1310-1343.








\bibitem{HERZ1971}
C. Herz,
{\it The theory of $p$-spaces with an application to convolution operators,
its second dual},
Trans. Amer. Math. Soc. \textbf{154} (1971) 69-82.













\bibitem{ILIE2013}
M. Ilie,
{\it A note on $p$-completely bounded homomorphisms of the Figà-Talamanca-Herz algebras},
J. Math. Anal. Appl. \textbf{419} (2014) 273-284.



\bibitem{ILIE2004}
M. Ilie,
{\it On Fourier algebra homomorphisms},
J. Funct. Anal. \textbf{213} (2004) 88-110.

\bibitem{ILIESPRONK2005}
 M. Ilie, N. Spronk,
{\it Completely bounded homomorphisms of the Fourier algebra}, 
J. Funct. Anal. \textbf{225} (2005) 480-499.   
   
   


\bibitem{ISTRATESCU1983}
V. I. Istratescu,
{\it Strict Convexity and Complex Strict Convexity: Theory and
Applications}, 
Taylor \& Francis Inc. 1983.
   
   






\bibitem{MIAO1996}
T. Miao,
{\it Compactness of a locally compact group $G$ and geometric properties of $A_p(G)$},
Canad. J. Math. \textbf{48} (1996) 127-1285.



\bibitem{OZTOPSPRONK2012} S. Oztop, N. Spronk,
{\it $p$-Operator space structure on Feichtinger-Fig\`a-Talamanca-Herz Segal algebra},
  J. Operator Theory \textbf{74} (2015), no. 1, 45-74. 



\bibitem{PIER1984} J.P. Pier,
{\it Amenable Locally Compact Groups}, Pure and Applied Math., 
 Wiley-Interscience, New York, 1984.













\bibitem{RUNDE2005}
 V. Runde,
{\it Representations of locally compact groups on ${QSL}_p$-spaces and a $p$-analog of the Fourier-Stieltjes algebra},
Pacific J. Math. \textbf{221} (2005) 379-397.






\bibitem{RUNDE2007}
V. Runde,
{\it Cohen-Host type idempotent theorems for representations on Banach spaces and applications to Fig\`a-Talamanca-Herz algebras},
J. Math. Anal. Appl. \textbf{329} (2007) 736-751.

\bibitem{shams}
M. Shams Yousefi,
{\it $p$-analog of the semigroup Fourier-Stieltjes algebras},
Iranian J. Math. Sci. and Inf. \textbf{10} (2) (2015) 55-66. 

	
	
	
	
	
	
	
\end{thebibliography}
\end{document}